\documentclass[final]{siamltex}
\usepackage{amsmath,amssymb,graphicx}
 \newtheorem{remark}[theorem]{Remark}

\title{Fourth order time-stepping for Kadomtsev-Petviashvili and 
Davey-Stewartson equations\thanks{We thank P.~Matthews, B.~Muite, 
A.~Ostermann, T. Schmelzer, who provided example codes, and
L.N.~Trefethen, who interested us in the subject,  
for helpful discussion and hints. This work has been supported by 
the project FroM-PDE funded by the European
Research Council through the Advanced Investigator Grant Scheme, 
the Conseil R\'egional de Bourgogne
via a FABER grant and the ANR via the program ANR-09-BLAN-0117-01.}}

\author{C.~Klein\thanks{Institut de Math\'ematiques de Bourgogne,
		Universit\'e de Bourgogne, 9 avenue Alain Savary, 21078 Dijon
		Cedex, France
    ({\tt christian.klein@u-bourgogne.fr})}
\and
K.~Roidot\thanks{Institut de Math\'ematiques de Bourgogne,
		Universit\'e de Bourgogne, 9 avenue Alain Savary, 21078 Dijon
		Cedex, France
    ({\tt kristelle.roidot@u-bourgogne.fr})}
}

\begin{document}
\maketitle

\begin{abstract}
    Purely dispersive partial differential 
       equations as the Korteweg-de Vries equation, the nonlinear 
	   Schr\"odinger equation and higher dimensional generalizations 
	   thereof	can have 
	   solutions 
	   which develop a zone of rapid modulated oscillations in 
	   the region where the corresponding dispersionless equations have 
	   shocks or blow-up. To numerically study such phenomena, fourth order time-stepping in combination 
	   with spectral methods is beneficial to 
	   resolve the steep gradients in the oscillatory 
	   region. We compare the performance of several fourth order 
	   methods for the Kadomtsev-Petviashvili  and the 
	   Davey-Stewartson equations, two integrable equations in 2+1 
	   dimensions:
	   exponential time-differencing, integrating factors, 
	   time-splitting, implicit Runge-Kutta and Driscoll's 
	   composite Runge-Kutta method. The 
	   accuracy in the numerical conservation of integrals of motion is discussed. 
\end{abstract}



\section{Introduction}
\begin{figure}[!htb]
\centering
\includegraphics[width=0.4\textwidth]{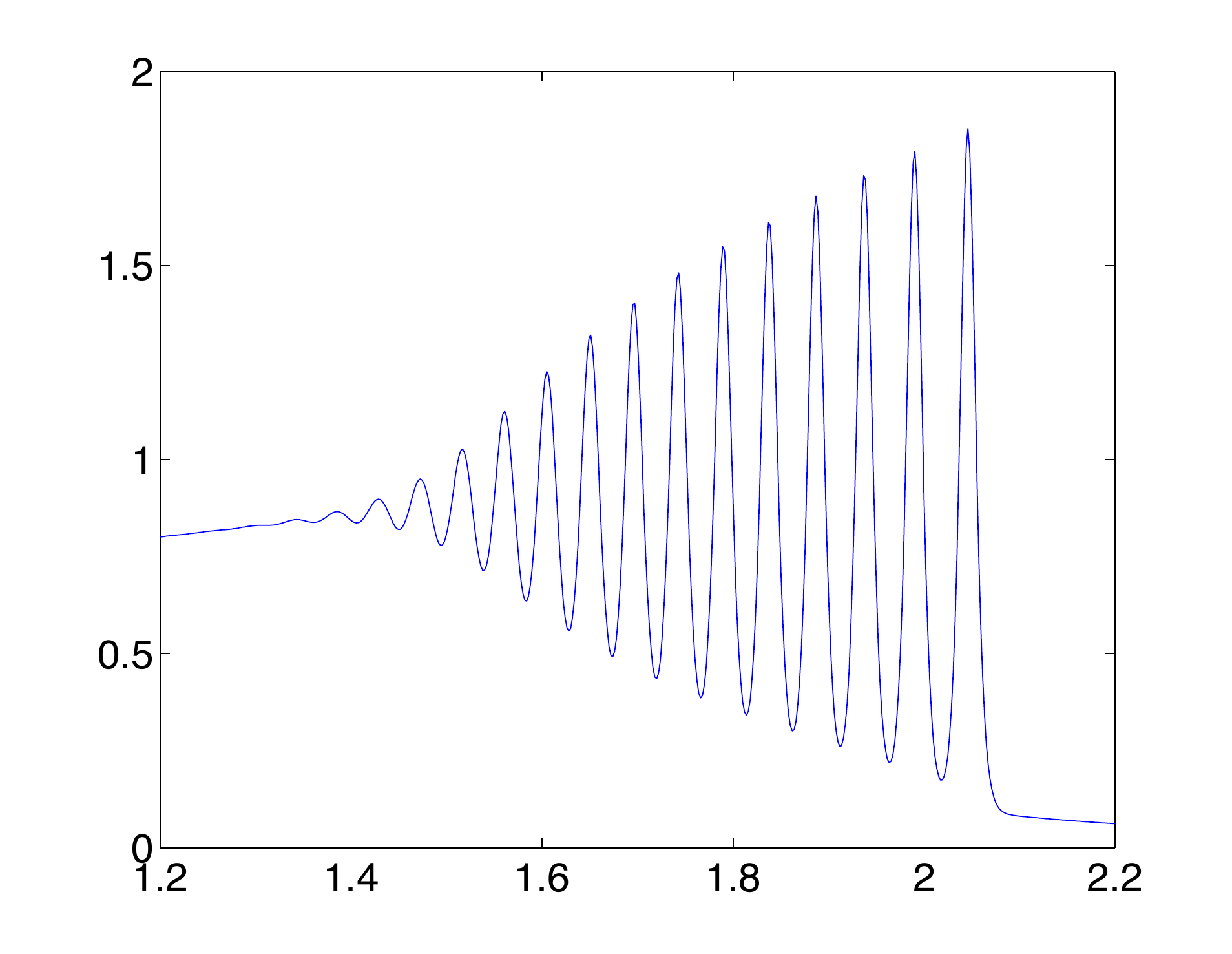}\quad
\includegraphics[width=0.4\textwidth]{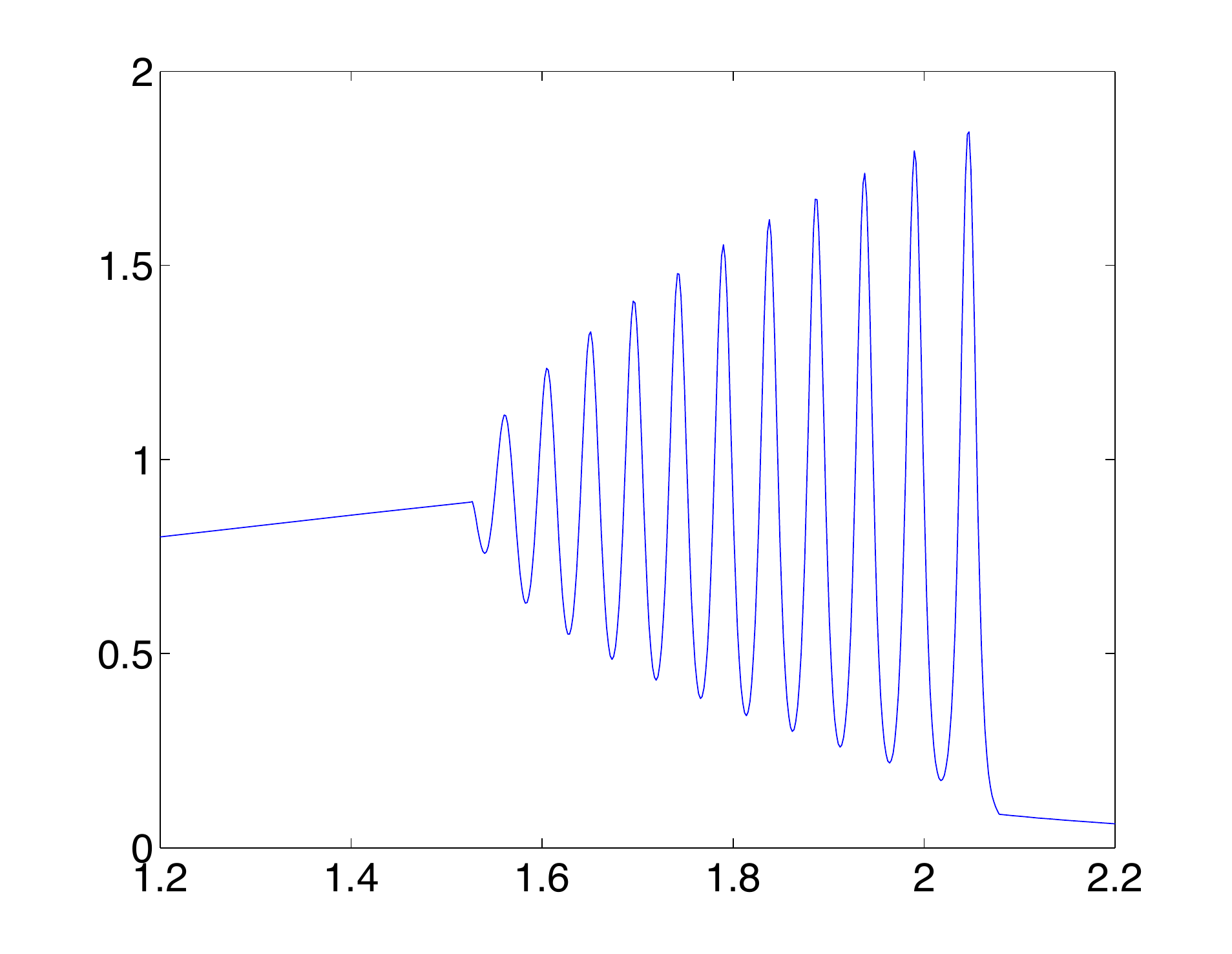}\newline
\caption{Numerical solution of KdV  for the initial data
 $u_0(x)=1/\cosh^2x$ (left) and the corresponding
 asymptotic solution (right) for  $t=0.35$ and  
 $\epsilon=10^{-2}$ (see \cite{GK}).}
\label{figposdata}
\end{figure}

Nonlinear dispersive partial differential equations (PDEs) play an important 
role in applications since they appear in many approximations to 
systems in hydrodynamics, nonlinear optics, acoustics, plasma 
physics, Bose-Einstein condensates among others. The most prominent members of the class are the 
celebrated Korteweg-de Vries (KdV) equation and the nonlinear 
Schr\"odinger (NLS) equation and higher dimensional generalizations 
of these. In addition to the importance of these equations in 
applications, there is also a considerable interest in the 
mathematical properties of their solutions. It is known that 
nonlinear dispersive PDEs without dissipation can have 
\emph{dispersive shock waves} \cite{GP}, i.e., regions of rapid modulated 
oscillations in the vicinity of shocks in the solutions to the corresponding 
dispersionless equations for the same initial data. Thus solutions to 
dispersive PDEs in general will not 
have a strong dispersionless limit as known from solutions to 
dissipative PDEs as 
the Burgers' equation in the limit of vanishing dissipation. An 
asymptotic description of these dispersive shocks is known for 
certain integrable PDEs as KdV \cite{LL,Ven,DVZ} and the NLS equation 
for certain classes of initial data 
\cite{JLM,KMM,TVZ}. For KdV an example is shown in 
Fig.~\ref{figposdata}, for details see \cite{GK}.

No such description is known for $2+1$-dimensional 
PDEs. In addition solutions to nonlinear dispersive PDEs can have 
\emph{blowup}, i.e., in finite time a loss of regularity of the solution with 
respect to the initial data. It is known for many of the 
PDEs under consideration when blowup can occur, but for the precise mechanism of the 
blowup often not even conjectures exist.

In view of the importance of the equations and the open mathematical 
questions, efficient numerical algorithms are needed to enable 
extensive numerical studies of the PDEs. The focus of the present work 
is to study a $2+1$-dimensional generalization of the KdV equation, the 
Kadomtsev-Petviashvili (KP) equation, and a $2+1$-dimensional 
generalization of the NLS equation, the Davey-Stewartson  (DS) equation. 
The former takes the form
\begin{equation}\label{e1}
\partial_{x}\left(\partial_{t}u+6u\partial_{x}u+\epsilon^{2}\partial_{xxx}u\right)+\lambda\partial_{yy}u=0,\,\,\lambda=\pm1
\end{equation}
where $(x,y,t)\in\mathbb{R}_{x}\times\mathbb{R}_{y}\times\mathbb{R}_{t}$
and where $\epsilon\ll1$ is a small scaling parameter. The limit 
$\epsilon\to0$ is the dispersionless limit. 
The case $\lambda=-1$ corresponds to the KP I model with a 
\emph{focusing} effect, and the case $\lambda=1$
corresponds to the KP II model with a \emph{defocusing} effect. The former is also known as 
the \emph{unstable} KP equation since the soliton solution of the KdV equation 
is lineary unstable for KP I, whereas it is lineary stable for the latter, which 
is therefore also known as the \emph{stable} KP equation. These 
stability issues were numerically studied in \cite{KS10}. It is 
an interesting result of the present paper that both KP equations 
have very similar numerical convergence properties despite completely 
different stability properties of their exact solutions in \cite{KS10}.
These equations appear in different fields of physics in
the study of
essentially one-dimensional wave phenomena with weak transverse effects, for example to model nonlinear dispersive
waves on the surface of fluids \cite{KP}. In this case, KP I is used when the surface tension is strong, and
KP II when the surface tension is weak. 
They also model sound waves in ferromagnetic media \cite{KPfer}, and nonlinear matter-wave pulses in
Bose-Einstein condensates \cite{KPbose}. The KP 
equation was introduced by Kadomtsev and Petviashvili in \cite{KP} to 
study the stability of the KdV soliton against weak transverse 
perturbations. It was shown to be completely integrable in 
\cite{Dry}. Higher dimensional generalizations of the KP equations, 
where the derivative $\partial_{yy}$ is replaced by the Laplacian 
in the transverse coordinates, 
$\Delta^{\perp}=\partial_{yy}+\partial_{zz}$, 
are important for instance in acoustics. The numerical problems to be 
expected there are the same as in the $2+1$-dimensional case studied 
here. 
  
The Davey-Stewartson system can be written in the form
\begin{equation}
    \label{DSII}
\begin{array}{ccc}
i\epsilon 
\partial_{t}u+\epsilon^{2}\partial_{xx}u-\alpha\epsilon^{2}\partial_{yy}u+2\rho\left(\Phi+\left|u\right|^{2}\right)u & = & 0,
\\
\partial_{xx}\Phi+\beta\partial_{yy}\Phi+2\partial_{xx}\left|u\right|^{2} & = & 0,
\end{array}
\end{equation}
where $\alpha$, $\beta$ and $\rho$ take the values $\pm1$, where
$\epsilon\ll1$ is again a small dispersion
parameter, and where $\Phi$ is a mean field. Since the $\epsilon$ has the 
same role as the $\hbar$ in the Schr\"odinger equation, the limit 
$\epsilon\to0$ is also called the semiclassical limit in this 
context. The DS equations are 
classified \cite{GS} according to the ellipticity or hyperbolicity of the 
operators in the first and second line. The case $\alpha=\beta$ is 
completely integrable  \cite{AH} and thus provides a $2+1$-dimensional 
generalization of the integrable NLS equation in $1+1$ dimensions. 
The integrable cases are
elliptic-hyperbolic  called DS I, and the hyperbolic-elliptic 
called DS II. For both there is a focusing ($\rho=-1$) and a 
defocusing ($\rho=1$) version. 
 In the following, we 
will only consider the case DS II ($\alpha=1$) since the mean field 
$\Phi$ is then obtained by inverting an elliptic operator. 
These DS systems model the evolution of weakly nonlinear water waves that
travel predominantly in one direction, but in which the wave amplitude
is modulated slowly in two horizontal directions \cite{DS}, \cite{DR}. They are also used
in plasma physics \cite{NAS,NAS94}, to describe the evolution of a plasma under the
action of a magnetic field. 

Since both KP and DS are completely integrable, there exist many
explicit solutions, which thus provide popular test cases for numerical
algorithms. But as we will show for the example of KP, these exact 
solutions, typically solitons, often test the equation in a 
regime where stiffness is not important. The main challenge in the study of critical phenomena as 
dispersive shocks and blowups is, however, the numerical resolution 
of strong gradients in the presence of which the above 
equations are very stiff. This implies that algorithms that perform well 
for solitons might not be efficient in the context studied here.

Since critical phenomena are generally believed to be independent of 
the chosen boundary conditions, we study a periodic 
setting\footnote{The boundary conditions will, however, in general 
influence convergence of the used numerical schemes. The 
restriction in this paper to periodic conditions is due to the studied 
problems.}. Such 
settings also include rapidly decreasing functions which can be 
periodically continued as smooth functions within the finite numerical precision. 
This allows one to approximate the spatial dependence 
via truncated Fourier series which leads for
the studied equations to large stiff systems of ODEs, see below. 
The use of 
Fourier methods not only gives spectral accuracy in the spatial 
coordinates, but also minimizes the introduction of numerical 
dissipation which is important in the study of dispersive effects.
In 
Fourier space, equations (\ref{e1}) and (\ref{DSII}) have the form
\begin{equation}
    v_{t}=\mathbf{L}v+\mathbf{N}(v,t)
    \label{utrans},
\end{equation}
where $v$ denotes the (discrete) Fourier transform of $u$, 
and where $\mathbf{L}$ and $\mathbf{N}$ denote linear and nonlinear 
operators, respectively. The resulting systems of ODEs are classical examples of stiff equations where the 
stiffness is related to the linear part $\mathbf{L}$ (it is 
a consequence of the distribution of the eigenvalues of 
$\mathbf{L}$), whereas the 
nonlinear part contains only low order derivatives. In the small dispersion 
limit, this stiffness is still present despite the small term 
$\epsilon^{2}$ in $\mathbf{L}$. This is due to the fact that the 
smaller $\epsilon$ is, the higher wavenumbers are needed to 
resolve the rapid oscillations. The first numerical studies of exact 
solutions to the KP equations were performed in \cite{WMGSS} and 
\cite{XS} and references therein. For the 
DS system similar studies were done in \cite{WW}. In \cite{BMS} 
blowup for DS was studied for the analytically known blowup solution 
by Ozawa \cite{Oza}.

There are several approaches to 
deal efficiently with equations of the form (\ref{utrans}) with a 
linear stiff part, implicit-explicit (IMEX), time splitting, 
integrating factor (IF), and deferred correction schemes as well as 
sliders and exponential time differencing. 
To avoid as much as possible a pollution of the Fourier coefficients 
by errors due to the finite difference schemes for the time 
integration and to allow the use of larger time steps, we mainly consider fourth order schemes.
While standard explicit schemes impose prohibitively small time 
steps due to stability requirements (for the studied examples the 
standard fourth order Runge-Kutta (RK) scheme did not converge for the 
used time steps), stable implicit schemes are in 
general 
computationally too expensive in $2+1$ dimensions. As an example of 
the latter we consider an implicit fourth order Runge-Kutta scheme. 
The focus of this paper is, however, to 
compare the performance of several explicit fourth order schemes mainly related to 
exponential integrators for various examples in a similar way as in 
the work by Kassam and Trefethen \cite{KassT} and in \cite{ckkdvnls} for 
KdV and NLS. 

The paper is organized as follows: In section 2 we briefly list the 
used numerical schemes, integrating factor methods, exponential time 
differencing, lineary implicit schemes, time splitting methods and 
implicit Runge-Kutta schemes. In 
section 3 we review some analytical facts for the KP equations and 
study  for each of the equations an exact solution and an example in the 
small  dispersion limit. In section 4 a similar 
analysis is presented for the semiclassical limit of the 
focusing and the defocusing DS II equation. The 
found numerical errors are compared in section 5 with the error 
indicated by a violation of the conservation of the $L_{2}$ norm by the numerical 
solution. In section 6 we add some concluding remarks 
and outline further directions of research. 
 
\section{Numerical Methods}
In this paper we are mainly interested in the numerical study of the KP and 
the DS II equations for Schwartzian initial data in the small 
dispersion limit. The latter implies that we can treat the problem as 
essentially periodic, and that we can use Fourier methods. After 
spatial discretization we thus face a system of ODEs of the form (\ref{utrans}). 
Since we need to resolve high wavenumbers, these systems will 
be in general rather large.  The PDEs studied here have high order 
derivatives in the linear part $\mathbf{L}$ of (\ref{utrans}), whereas the 
nonlinear part $\mathbf{N}$ contains only first derivatives. This means that 
the stiffness in these systems is due to the linear part. The latter 
will thus be treated with adapted methods detailed below, whereas 
standard methods can be used for the nonlinear part. We restrict the analysis 
to moderate values of the dispersion parameter	to be able to 
study the dependence of the different 
schemes on the time step in finite CPU time. For smaller values of $\epsilon$ see for 
instance \cite{KSM}.

We will compare several numerical schemes for equations of the form 
(\ref{utrans}) as in \cite{KassT} and \cite{ckkdvnls}. The PDEs are 
studied for a periodic setting with periods $2\pi L_{x}$ and $2\pi 
L_{y}$ in $x$ and $y$ respectively. We give the numerical error in 
dependence of the number $N_{t}$ of time steps as well as the actual CPU time as measured by MATLAB 
(all computations are done on a machine with Intel  `Nehalem' 
processors with 2.93 GHz with codes in MATLAB 7.10). The goal 
is to provide some indication on the actual performance of the codes 
in practical applications. 
Since MATLAB is using in general a mixture of interpreted and precompiled embedded code, 
a comparison of computing times is not unproblematic. However, it can 
be done in the present context since the main computational cost  is 
due to 
two-dimensional fast Fourier transformations (FFT). For the KP equations all considered approaches 
(with the exception of the Hochbruck-Ostermann ETD scheme which 
uses 8 FFT commands per time step) use 6 (embedded) FFT commands per time step as 
was already pointed out in \cite{D}. For the DS II equation, these 
numbers are doubled since the computation of the mean field $\Phi$ 
takes another FFT/IFFT pair per intermediate step. Note that an 
additional FFT/IFFT pair per time step is needed in both cases to switch 
between physical and Fourier space (we are interested in a solution 
in physical space, but the schemes are formulated for the Fourier 
transforms). 
The $\phi$-functions in the ETD schemes are
also computed via FFT. It can be seen that this 
can be done with machine precision in a very efficient way. 
Since the $\phi$-functions have to be obtained only once in the computation and since the 
studied problems are computationally demanding, this only has a negligible 
effect on the total CPU time in the experiments.
The numerical error is the $L_{2}$ norm of the difference of 
the numerical solution and an exact or reference solution, normalized by the	
$L_{2}$ norm of the initial data. It is denoted by $\Delta_{2}$.

\subsection{Integrating Factor Methods (IF)}

These methods appeared first in the work of Lawson 
\cite{law}, see \cite{MW} for a review.
He suggested to take care of the stiff linear part of 
equation (\ref{utrans}) by
using a change of the dependent variables (also called the Lawson transformation)
$w(t)=e^{-\mathbf{L}t}v(t)$.
Equation (\ref{utrans}) becomes
\begin{equation}\label{e6}
w'(t)=e^{-\mathbf{L}t}\mathbf{N}(v,t)
\end{equation}
for which we use a fourth order Runge-Kutta (RK) scheme. Hochbruck and 
Ostermann \cite{HO} showed that this IF method has classical 
order four, but not what they call \emph{stiff order} four. Loosely 
speaking there can be additional contributions to the error of much 
lower order in the time step due to the stiffness, i.e., in the case 
of large norm $||L||$ in the integrating factor. They could show that 
the scheme used here can reduce to first order for semilinear parabolic 
problems. The Hochbruck-Ostermann approach uses 
semigroups. An extension of their theory to 
hyperbolic equations is possible via $C_{0}$ semigroups, which can 
lead, however, to a slightly lower order of convergence. But the results in 
\cite{BS}, \cite{ckkdvnls} indicate that similar convergence rates are to be 
expected for hyperbolic equations (both KP and DS are hyperbolic in 
the sense that the matrix $\mathbf{L}$ appearing in (\ref{utrans}) has purely 
imaginary eigenvalues).

\subsection{Driscoll's composite Runge-Kutta Method}
The idea of IMEX methods (see e.g.~\cite{CK} for KdV) is the use of a 
stable implicit method for the linear part of the equation 
(\ref{utrans}) and an explicit scheme for the nonlinear part which is 
assumed to be non-stiff. In \cite{KassT} such schemes did not 
perform satisfactorily for dispersive PDEs which is why we only 
consider a more sophisticated variant here. Fornberg and Driscoll 
\cite{Forn} provided an interesting generalization of IMEX by splitting 
also the linear part of the equation in Fourier space into regimes of high, 
medium, and low wavenumbers, and by using adapted numerical schemes in 
each of them. They considered the NLS equation as an example. 
Driscoll's \cite{D} idea  was to split the linear part of the equation in
Fourier space just into regimes of high and low wavenumbers. He
used the fourth order RK integrator for the low wavenumbers and the
lineary implicit RK method of order three for the high wavenumbers.
He showed that this method is in practice of fourth order over a wide range of step
sizes. We confirm this here for the cases where the method converges, 
which it fails to do, however, sometimes for very stiff problems.
In
particular, he used this method for the KP II equation at the two 
phase solution we will also discuss in this paper as a test case. We 
call the method DCRK in the following.

\subsection{Exponential Time Differencing Methods}

Exponential time differencing schemes were developed originally by 
Certaine \cite{Cer} in the 60s,  see 
\cite{MW} and \cite{HO09} for comprehensive 
reviews of ETD methods and their  history. The basic idea is to use equidistant time 
steps $h$ and to integrate equation (\ref{utrans}) exactly between 
the time steps $t_{n}$ and $t_{n+1}$ with respect to $t$. With $v(t_{n}) = v_{n}$ and 
$v(t_{n+1})=v_{n+1}$, we get 
$$
v_{n+1}=e^{\mathbf{L}h}v_{n}+\int_{0}^{h} e^{\mathbf{L}(h-\tau)}
    \mathbf{N}(v(t_{n}+\tau),t_{n}+\tau)d\tau.
$$
The integral will be computed in an approximate way for which 
different schemes exist. We use here only  
Runge-Kutta schemes of classical order 4,  Cox-Matthews 
\cite{CM}, Krogstad \cite{K} and  Hochbruck-Ostermann  \cite{HO}. 
The latter showed that the stiff order of 
the Cox-Matthews scheme is only two, and the one of Krogstad's is 
three. Notice that both schemes can be of stiff order four if the 
studied system satisfies a certain number of non-trivial auxiliary conditions, see 
\cite{HO}. As the numerical 
tests show in the following, order reduction can be observed in some 
cases, but not in all. We will speak in the following of the 
\emph{stiff regime} of an equation where order reduction can be 
observed, or where certain schemes do not converge, and of the 
\emph{non-stiff} regime, where this is not the case. Note that this 
is not a rigorous definition. Both Cox-Matthews' and Krogstad's schemes
are, however, four-stage methods, whereas the 
Hochbruck-Ostermann method  is a five-stage method that has stiff order 
four. Thus all these methods should show the same convergence rate in 
the non-stiff regime, but could differ for some problems in the stiff regime. Notice that 
these results \cite{HO} were established for parabolic PDEs, and that 
the applicability for hyperbolic PDEs of the type studied here is 
not obvious. One of the purposes of our study is to get some 
experimental insight whether the Hochbruck-Ostermann theory holds 
also in this case.

The main 
technical problem in the use of  ETD schemes is the efficient and accurate 
numerical evaluation of the functions 
$$    \phi_{i}(z) = 
    \frac{1}{(i-1)!}\int_{0}^{1}e^{(1-\tau)z}\tau^{i-1}d\tau, \quad 
    i=1,2,3,4,$$
i.e., functions of the form $(e^{z}-1)/z$ and higher order 
generalizations thereof, where one has to avoid 
cancellation errors. Kassam and Trefethen \cite{KassT} used complex 
contour integrals to compute these functions. The approach is 
straight forward for diagonal operators $\mathbf{L}$ that occur here 
because of the use of Fourier methods: 
one considers a unit circle around each point $z$ and computes the 
contour integral with the trapezoidal rule which is known to be a 
spectral method in this case. 
Schmelzer \cite{schme} made this approach more efficient by 
using the complex contour approach only for values of $z$ close to 
the pole, e.g\ with 
$|z|<1/2$. For the same values of $z$ the functions $\phi_{i}$ can be 
computed via a Taylor series. These two independent and very 
efficient approaches 
allow a control of the accuracy. We find that just 16 Fourier 
modes in the computation of the 
complex contour integral are sufficient to determine the functions 
$\phi_{i}$ to the order of machine precision. Thus we avoid problems 
reported in \cite{BS}, where machine precision could not be reached 
by ETD schemes due to inaccuracies in the determination of the $\phi$-functions. 
The computation of these functions takes 
only negligible time for the $2+1$-dimensional equations studied here, 
especially since it has to be done only once 
during the time evolution. We find that ETD as implemented in this 
way has the same computational costs as the other used schemes.

\subsection{Splitting Methods}

Splitting methods are convenient if an equation can be split into two or
more equations which can be directly integrated.
The  motivation for these methods is the Trotter-Kato formula
\cite{TK,Ka}
\begin{equation}\label{e11}
 {\lim}_{n\rightarrow\infty}\left(e^{-tA/n}e^{-tB/n}\right)^{n}=e^{-t\left(A+B\right)}
\end{equation}
where $A$ and $B$ are certain unbounded linear operators, for details 
see \cite{Ka}.
In particular this includes the cases studied by Bagrinovskii and
Godunov in \cite{BG} and by Strang \cite{ST}. For hyperbolic equations,
first references are Tappert \cite{Tap} and Hardin and Tappert \cite{HT}
who introduced the split step method for the NLS
equation.

The idea of these methods for an equation of the form $u_{t}=\left(A+B\right)u$ is to write the solution in the
form
\[
u(t)=\exp(c_{1}tA)\exp(d_{1}tB)\exp(c_{2}tA)\exp(d_{2}tB)\cdots\exp(c_{k}tA)\exp(d_{k}tB)u(0)
\]
where $(c_{1},\,\ldots,\, c_{k})$ and $(d_{1},\,\ldots,\, d_{k})$
are sets of real numbers that represent fractional time steps.  
Yoshida \cite{Y} gave an approach 
which produces split step methods of any even order.

The KP equation can be split into
\begin{eqnarray}
u_{t}+6uu_{x}=0,  \label{kpho}
\\
(\mathcal{F}[u])_{t}-ik_{x}^{3}\mathcal{F}[u]+\lambda 
\frac{ik_{y}^{2}}{k_{x}}\mathcal{F}[u]=0, 
\label{kpsplit}              
\end{eqnarray}
where 
here and in the following we write the 2-dimensional Fourier 
transform of $u$ in the form
\begin{equation}
    \mathcal{F}[u]:=\int_{\mathbb{R}^{2}}^{}u(x,y,t)e^{-ik_{x}x-ik_{y}y}dxdy.
    \label{fourier2d}
\end{equation}
The Hopf equation (\ref{kpho}) can be integrated in implicit form with the method of 
characteristics, and the linear equation in Fourier space (\ref{kpsplit})
can be directly integrated, but the implicit form of the solution of the 
former makes an iteration with interpolation to the characteristic 
coordinates necessary that is computationally too 
expensive. Therefore we consider splitting here only for the DS 
equation. The latter can be split into
\begin{eqnarray}
i\epsilon u_{t}=\epsilon^{2}(-u_{xx}+\alpha u_{yy}),\,\,\, 
\Phi_{xx}+\alpha\Phi_{yy}+2\left(\left|u\right|^{2}\right)_{xx}=0 ,
\\
i\epsilon u_{t}= -2\rho\left(\Phi+\left|u\right|^{2}\right)u,  
\label{2.4}              
\end{eqnarray}
which are explicitly integrable, the first two in Fourier space, 
equation (\ref{2.4}) in physical space since 
$|u|^{2}$ is a constant in time for this equation. Convergence of second order 
splitting along these lines was studied in \cite{BMS}. We study here 
second and fourth order splitting schemes for DS as given in \cite{Y}. 

\subsection{Implicit Runge Kutta Scheme}

The general formulation of an $s$-stage Runge--Kutta method for the initial value problem
$y'=f(y,t),\,\,\,\,y(t_0)=y_0$ is the following:
\begin{eqnarray}
 y_{n+1} = y_{n} + h      \underset{i=1}{\overset{s}{\sum}} \, 
 b_{i}K_{i}, \\
 K_{i} = f\left(t_{n}+c_ {i}h,\,y_{n}+h  
 \underset{j=1}{\overset{s}{\sum}} \, a_{ij}K_{j}\right),
\end{eqnarray}
where $b_i,\,a_{ij},\,\,i,j=1,...,s$ are real numbers and
$c_i=   \underset{j=1}{\overset{s}{\sum}} \, a_{ij}$.

For the implicit Runge--Kutta scheme of order 4 
(IRK4) used here (Hammer-Hollingsworth method), we have 
$c_{1}=\frac{1}{2}-\frac{\sqrt{3}}{6}$, 
$c_{2}=\frac{1}{2}+\frac{\sqrt{3}}{6}$, $a_{11}=a_{22}=1/4$,
$a_{12}=\frac{1}{4}-\frac{\sqrt{3}}{6}$, 
$a_{21}=\frac{1}{4}+\frac{\sqrt{3}}{6}$ and $b_{1}=b_{2}=1/2$. This 
scheme is also known as the 2-stage Gauss method. It is of classical 
order 4, but stage order 2. This is the reason why an order reduction 
to second order can be observed in certain examples. 

The implicit character of this method requires the iterative solution of a high dimensional system 
at every step which is done via a simplified Newton method. 
For the studied examples in the form (\ref{utrans}), we have to solve 
equations of the form 
$$y = \mathbf{A}y+b(y)$$
for $y$, where $\mathbf{A}$ is a linear operator independent of $y$, and where 
$b$ is a vector with a nonlinear dependence on $y$. These are solved 
iteratively in the form
$$y_{n+1} = (1-\mathbf{A})^{-1}b(y_{n}).$$
By treating the linear part that is responsible for the stiffness 
explicitly as in an IMEX scheme, the iteration converges in general 
quickly. Without taking explicit care of the linear part, convergence 
will be extremely slow.
The iteration 
is stopped once the $L_{\infty}$ norm of the difference between 
consecutive iterates is smaller than some threshold (in practice we 
work with a threshold of $10^{-8}$). Per iteration the computational 
cost is essentially 2 FFT/IFFT pairs. Thus the IRK4 scheme can be 
competitive with the above explicit methods which take 3 or 4 
FFT/IFFT pairs per time step if not more than 2-3 iterations are 
needed per time step. This can happen in the below examples in the 
non-stiff regime, but is not the case in the stiff regime. We only 
test this scheme where its inclusion appears interesting and where it is 
computationally not too expensive.


%
%

\section{Kadomtsev-Petviashvili Equation}
In this section we study the efficiency of the above mentioned numerical 
schemes in solving  Cauchy problems for the KP equations. We first 
review some analytic facts about KP I and KP II which are important 
in this context. Since the KP equations are completely integrable, 
exact solutions exist that can be used as test cases for the codes. 
We compare the performance of the codes for the exact solutions and a 
typical example in the small dispersion limit.

\subsection{Analytic Properties of the KP Equations}
We will collect here some analytic aspects of the KP equations which 
will be important for an understanding of several issues in the 
numerical solution of Cauchy problems for the KP equations, see 
\cite{KS10} for a recent review and references therein. 

In this paper we will look for KP solutions that are periodic in $x$ 
and $y$, i.e., for solutions on $\mathbb{T}^{2}\times \mathbb{R}$. 
This includes for numerical purposes the case of rapidly decreasing 
functions in the Schwartz space $\mathcal{S}(\mathbb{R}^{2})$ if the 
periods are chosen large enough that $|u|$ is smaller than machine 
precision (we work with double precision throughout the article) at 
the boundaries of the computational domain. 
Notice, however, that solutions to Cauchy problems with Schwartzian 
initial data $u_{0}(x,y)$ will not stay in 
$\mathcal{S}(\mathbb{R}^{2})$ unless $u_{0}(x,y)$ satisfies an 
infinite number of constraints. This behaviour can be already seen on 
the level of the linearized KP equation, see e.g.~\cite{BPP,KSM}, where 
the Green's function implies a slow algebraic decrease in $y$ towards 
infinity. This leads to the formation of \emph{tails} with an 
algebraic decrease to infinity for generic Schwartzian initial data. 
The amplitude of these effects grows with time (see for instance \cite{KSM}).
In our periodic setting this will give rise to \emph{echoes} and a weak Gibbs phenomenon at 
the boundaries of the computational domain. The latter implies that 
we cannot easily reach machine precision as in the KdV case unless we 
use considerably larger domains. As can be seen from computations in 
the small dispersion limit below and the Fourier coefficients in 
sect.~5, we can nonetheless reach an 
accuracy of better than $10^{-10}$ on the chosen domain. For higher 
precisions  and larger values of $t$,  the 
Gibbs phenomena due to the algebraic tails become important.

The KP equation (\ref{e1}) is not in the standard form 
for a Cauchy problem. As discussed in \cite{KSM} $t$ is not a 
timelike but characteristic coordinate if the dispersionless KP 
equation ($\epsilon=0$) is considered as a standard second order PDE. 
In practice one is, however, interested in the Cauchy problem for 
$t$. To this end one writes (\ref{e1}) in \emph{evolutionary 
form} 
\begin{equation}
   \partial_{t}u+6u\partial_{x}u+\epsilon^{2}\partial_{xxx}u=
   -\lambda\partial_{x}^{-1}
   \partial_{yy}u,\,\,\lambda=\pm1
    \label{kpevol}.
\end{equation}
Equations (\ref{kpevol}) and (\ref{e1}) are equivalent for certain 
classes of boundary conditions as periodic or rapidly decreasing at 
infinity. Since we will always impose periodic boundary conditions in 
the following, both forms of the KP equation are equivalent for our 
purposes.
The antiderivative $\partial_{x}^{-1}$ is to be understood as 
the Fourier multiplier with the singular symbol $-i/k_{x}$. In the 
numerical computation this multiplier is regularized in standard way 
(similar to the Dirac regularization of $1/x$ by adding an arbitrary 
small imaginary part $i0$ to $x$) as 
$$\frac{-i}{k_{x}+i\lambda\delta},$$
where we choose $\delta=\mbox{eps}=2^{-52}\sim 2.2*10^{-16}$. Since 
the typical precision to be achieved in the studied examples with double precision is of the 
order $10^{-14}$ because of rounding errors, this is essentially 
equivalent to adding a numerical zero (see also the discussion in 
\cite{KSM}). 

The divergence structure of the KP equations has the consequence that 
\begin{equation}
    \int_{\mathbb{T}}^{}\partial_{yy}u(x,y,t)dx=0,\quad \forall t>0
    \label{const},
\end{equation}
even if this constraint is not satisfied for the initial data 
$u_{0}(x,y)$. It was shown in \cite{FS,MST} that the solution to a 
Cauchy problem not satisfying the constraint will not be smooth in 
time for $t=0$. Numerical experiments in \cite{KSM} indicate that 
the solution after an arbitrary small time step will develop an 
infinite `trench' the integral over which just ensures that 
(\ref{const}) is fulfilled. To propagate such initial data, the above 
regularization is in fact needed (for data satisfying the constraint 
via the condition that the Fourier coefficients for $k_{x}=0$ vanish, this 
property could be just imposed at each time step). 
The infinite trench due to initial data not satisfying the 
constraint implies a rather strong Gibbs 
phenomenon. To avoid the related problems, we always consider initial 
data that satisfy (\ref{const}). A possible way to achieve this is to 
consider data that are $x$-derivatives of periodic or Schwartzian 
functions. 

The complete integrability of the KP equations implies that efficient 
tools exist for the generation of exact solutions. We always put 
$\epsilon=1$ for the exact solutions.  For a recent 
review of the integrable aspects of KP see \cite{Fok09}. The most 
popular KP solutions are \emph{line solitons}, i.e., solutions 
localized in one spatial direction and infinitely extended in 
another. Such solutions typically have an angle not equal to 0 or 90 
degrees with the boundaries of the computational domain, which leads to strong Gibbs 
phenomena. This implies that these solutions are not a good test case 
for a periodic setting. If the angle is 0 or 90 degrees, the solution only depends 
on one of the spatial variables and thus does not test a true 2d 
code. There exists a \emph{lump} soliton for KP I which is localized in all 
spatial directions, but only with algebraic fall off. This would again lead to 
strong Gibbs phenomena in our setting.

However a solution due to Zaitsev \cite{Zai} to the KP I equation is localized 
in one direction and periodic in the second (a transformation of the 
form $x\to i x$, $y\to iy$ exchanges these two directions). 
It has the form
\begin{equation}
    u(\xi,y)=2\alpha^{2}\frac{1-\beta\cosh(\alpha\xi)\cos(\delta y)}{\left(\cosh(\alpha\xi)-\beta\cos(\delta y)\right)^{2}}
    \label{zait}
\end{equation}
where 
\[
\xi=x-ct , \, \, c=\alpha^{2}\frac{4-\beta^{2}}{1-\beta^{2}}, \, \, \, \, \mbox{and} \, \, \, \, 
\delta=\sqrt{\frac{3}{1-\beta^{2}}}\alpha^{2}.
\]
This solution  is localized in $x$, 
periodic in $y$, and unstable as discussed in \cite{KS10}.

Algebro-geometric solutions to the KP equation can be constructed  on an
arbitrary compact Riemann surface, see e.g.~\cite{Dub}, 
\cite{FK}. These solutions are in general almost 
periodic. Solutions on genus 2 surfaces, which are all hyperelliptic, 
are exactly periodic, but in general not in both $x$ and $y$. A 
doubly periodic solution with exactly this property of KP II of genus 2 can be written as
\begin{equation}
u(x,\, y,\, t)=2\frac{\partial^{2}}{\partial x^{2}}\ln\theta\left(\varphi_{1},\,\varphi_{2};\, B\right)
\end{equation}
where $\theta\left(\varphi_{1},\,\varphi_{2};\, B\right)$ is defined by the double Fourier series
\begin{equation}
\theta\left(\varphi_{1},\,\varphi_{2};\, B\right)=\overset{\infty}{\underset{m_{1}=-\infty}
{\sum}}\overset{\infty}{\underset{m_{2}=-\infty}{\sum}}e^{\frac{1}{2}m^{T}Bm+im^{T}\varphi}
\end{equation}
where $m^{T}=\left(m_{1},\, m_{2}\right)$, and where $B$ is a $2\times2$
symmetric, negative-definite Riemann matrix 
\[
B=\left(\begin{array}{cc}
b & b\lambda\\
b\lambda & b\lambda^{2}+d\end{array}\right),\,\,\,\,
\mbox{with real parameters}\,\,\,\lambda\neq0,\, b\,\,\,\mbox{and}\,\,\, d.
\]
 The phase variable $\varphi$ has the form $\varphi_{j}=\mu_{j}x+\nu_{j}y+\omega_{j}t+\varphi_{j,0},\,\,\, j=1,2.$
The solution travels as the 
Zaitsev solution with constant speed in $x$-direction.

\begin{remark}
    The standard 4th order Runge-Kutta scheme did not converge for 
    any of the studied examples for the used time steps. The reason is 
    that the Fourier multiplier $-i/k_{x}$ imposes very strong 
    stability restrictions on the scheme. 
\end{remark}

\subsection{Numerical solution of Cauchy problems for the KP I equation}

\paragraph{Zaitsev solution}

\begin{figure}[htb!]
 \centering
 \includegraphics[width=0.49\textwidth]{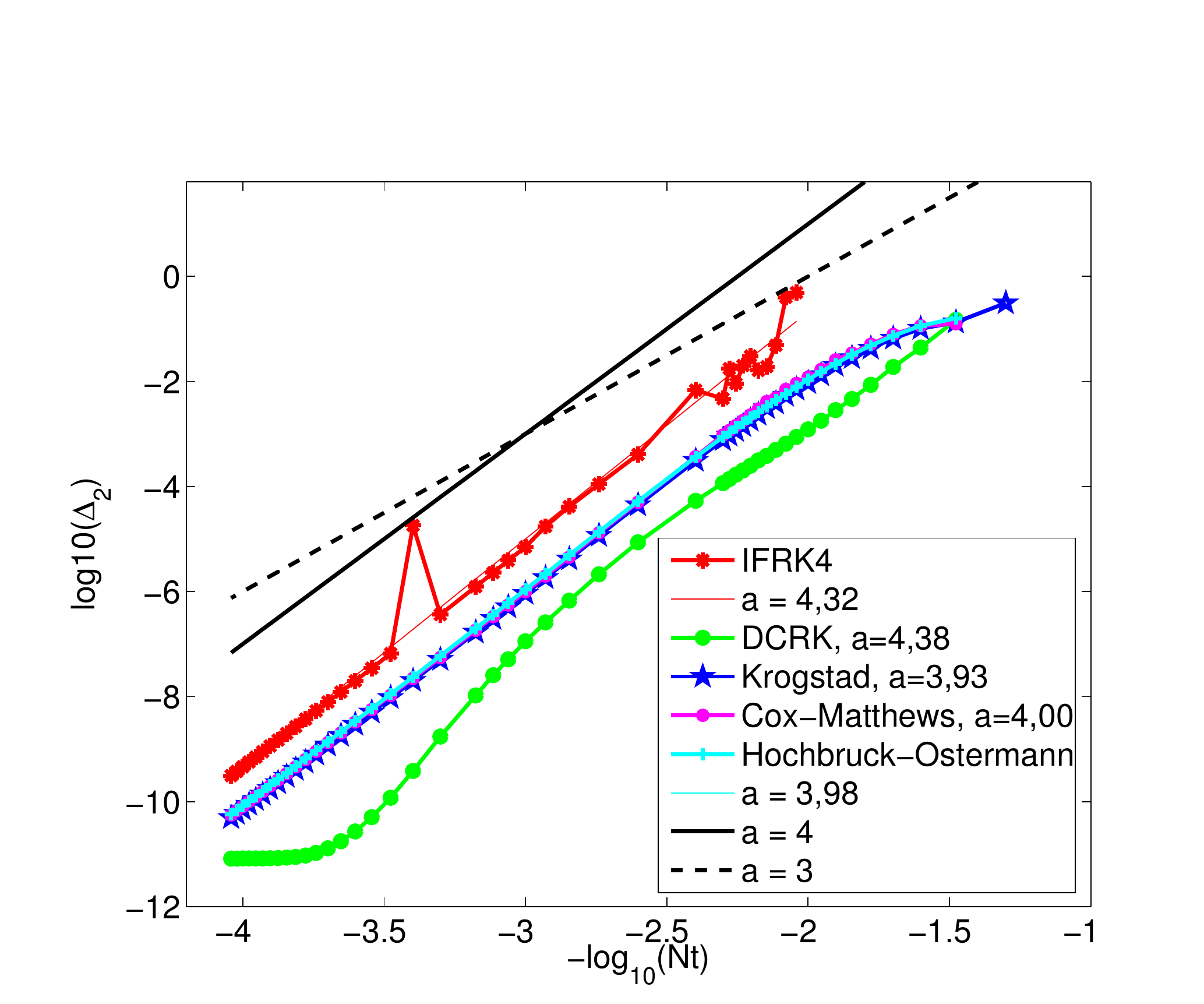}
 \includegraphics[width=0.49\textwidth]{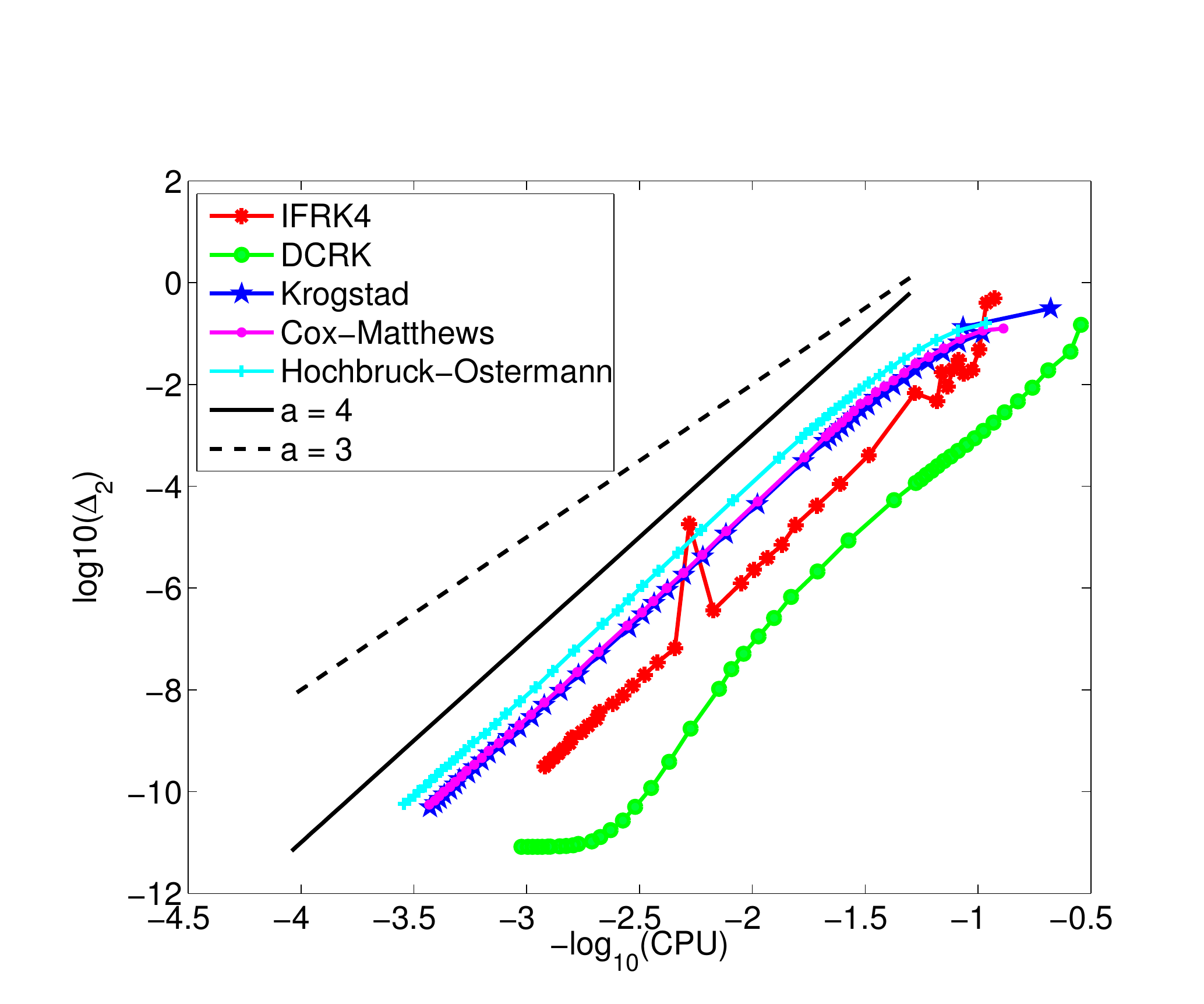}
 \caption{Normalized $L_{2}$ norm of the numerical error 
 in the solution to the KP I equation with initial data given by the 
 Zaitsev solution for several numerical
 methods, as a function of $N_{t}$ (left) and as a function of CPU 
 time (right).}
 \label{figzaitreg}
\end{figure}

We first study the case of the Zaitsev solution (\ref{zait}) with 
$\alpha=1$ and $\beta=0.5$. Notice that this 
solution is unstable against small perturbations
as shown numerically in \cite{KS10}, but that it 
can be propagated with the expected numerical precision by the used 
codes. As initial data we 
take the solution centered at $-L_{x}/2$ (we use $L_{x}=L_{y}=5$) 
and propagate it until it 
reaches $L_{x}/2$. 
The computation is carried out with $2^{11}\times2^{9}$ points for
$(x, y)\in[-5\pi,\,5\pi]\times[-5\pi,\,5\pi]$ and $t\leq1$.
The decrease of the numerical error is shown in 
Fig.~\ref{figzaitreg} in dependence of the time step and 
in dependence of the CPU time.
A linear regression analysis in a double logarithmic plot 
($\log_{10}\Delta_{2}=-a\log_{10}N_{t}+b$)
is presented in Fig.~\ref{figzaitreg}, where we can see that all schemes 
show a fourth order behavior:   we find $a=4.32$ for the Integrating Factor method, $a=4.38$ for
DCRK method, $a=3.93$ for Krogstad's ETD scheme, $a=4$ for 
the Cox-Matthews scheme, and $a=3.98$ for the Hochbruck-Ostermann scheme.
In this context the DCRK method performs best, 
followed by the ETD schemes that have almost 
identical performance (though the Hochbruck-Ostermann method uses 
more internal stages and thus more CPU time in 
Fig.~\ref{figzaitreg}). It can also be seen 
that the various schemes do not show the phenomenon of order 
reduction as discussed in \cite{HO}, which implies that the Zaitsev 
solution tests the codes in a non-stiff regime of the KP I equation. 

\begin{figure}[htb!]
 \centering
 \includegraphics[width=\textwidth]{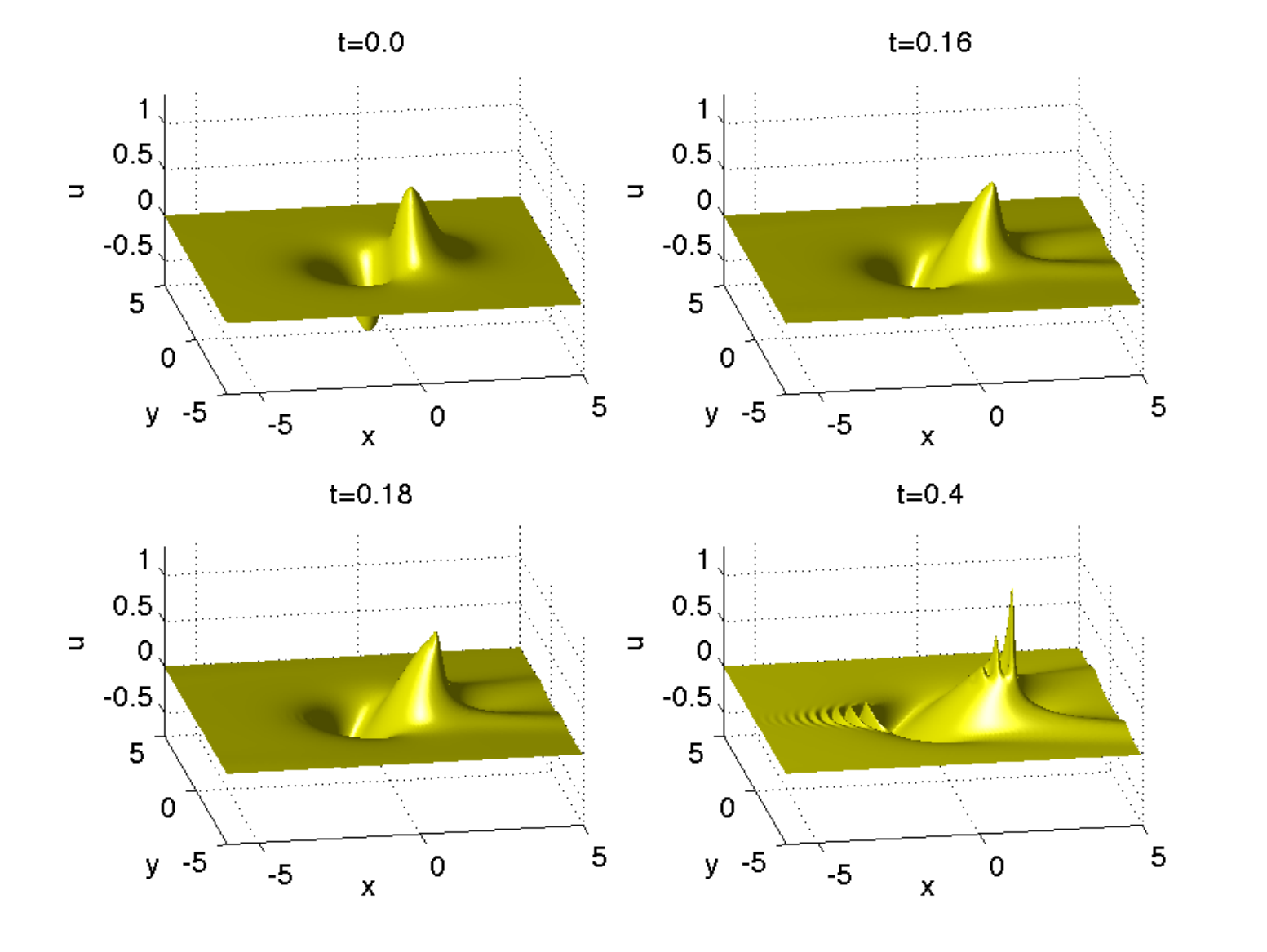}
 \caption{Solution to the KP I equation for the initial data $u_{0}=-\partial_{x}\mbox{sech}^{2}(R)\,\,\,\,\,\,
\mbox{where}\,\,\, R=\sqrt{x^{2}+y^{2}}$ for several values of $t$.}
 \label{figKPI}
\end{figure}

\paragraph{Small dispersion limit for KP I}
To study KP solutions in the limit of small dispersion 
($\epsilon\to0$), we consider Schwartzian initial data satisfying the constraint 
(\ref{const}).
As in \cite{KSM} we consider 
data of the form 
\begin{equation}\label{e12}
u_{0}(x,\, y)=-\partial_{x}\mbox{sech}^{2}(R)\,\,\,\,\,\,
\mbox{where}\,\,\, R=\sqrt{x^{2}+y^{2}}.
\end{equation}
By numerically solving the dispersionless KP equation (put 
$\epsilon=0$ in (\ref{e1})), we determine the critical time of the 
appearance of a gradient catastrophe by the breaking of the code, see 
\cite{KSM}. To study dispersive shocks, we run the KP codes for some 
time larger than this critical time. The solution can be seen 
in Fig.~\ref{figKPI}. It develops tails with algebraic fall 
off towards infinity. The wave fronts steepen on both sides of the 
origin. In the regions of strong gradients, rapid modulated 
oscillations appear. For a detailed discussion, see \cite{KSM}.

\begin{figure}[htb!]
 \centering
 \includegraphics[width=0.49\textwidth]{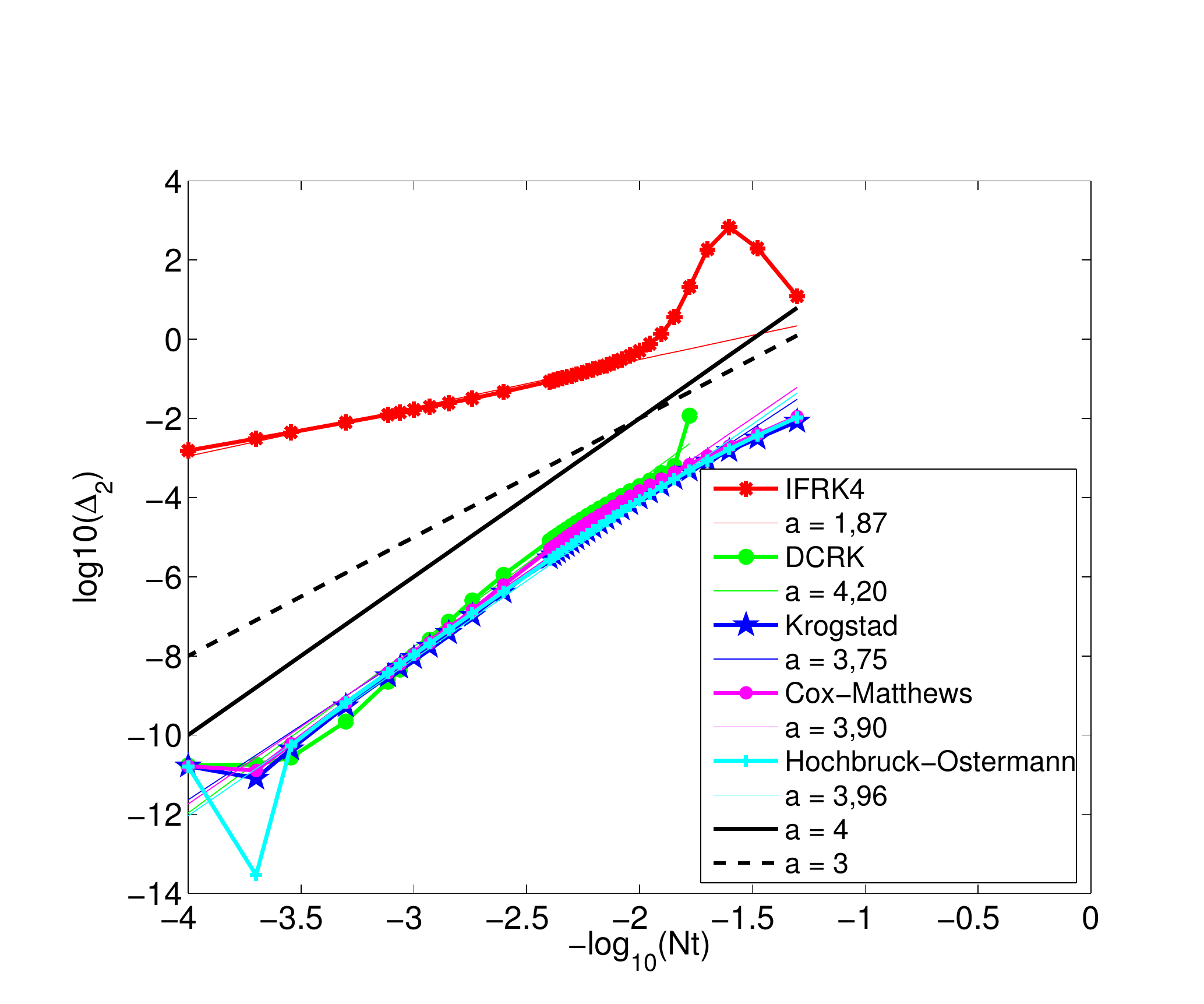}
 \includegraphics[width=0.49\textwidth]{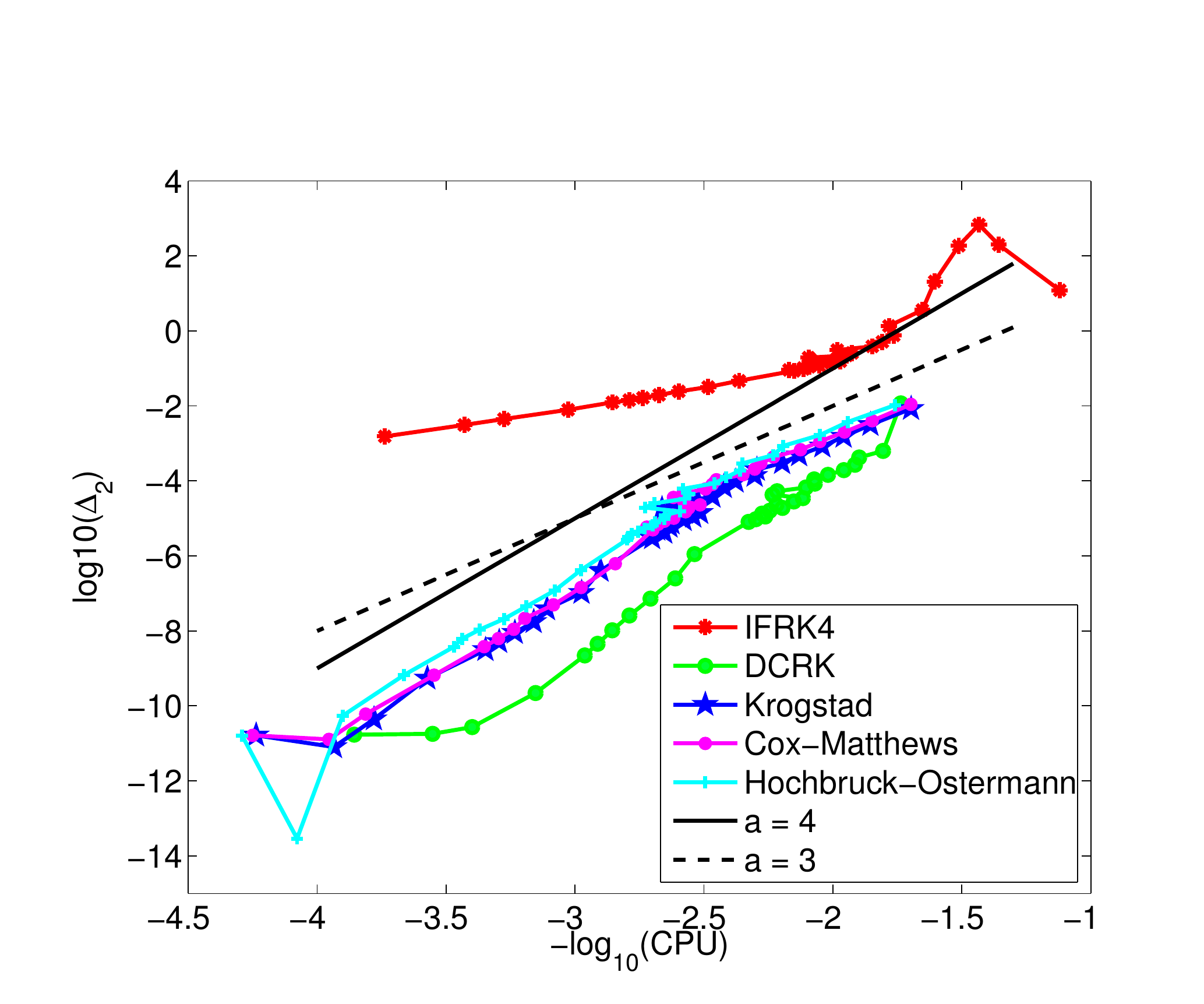}
 \caption{Normalized $L_{2}$ norm of the numerical error for the 
  solution shown in Fig.~\ref{figKPI} for several numerical
  methods, as a function of $N_{t}$  (left) and as a function of CPU 
 time (right).}
 \label{figKPIsmallreg}
\end{figure}

The computation
is carried out with $2^{11}\times2^{9}$ points for
$(x, y)\in[-5\pi,\,5\pi]\times[-5\pi,\,5\pi]$,
$\epsilon=0.1$ and $t\leq0.4$. As a reference solution, we consider
the solution calculated with the Hochbruck-Ostermann method with $N_{t}=5000$
time steps. The normalized $L_{2}$ norm of the
difference between this reference solution and the numerical solution
is shown in Fig.~\ref{figKPIsmallreg} in dependence on the time step 
with a regression analysis and  in dependence on the CPU time.
Here we can see clearly the phenomenon of 
order reduction established analytically for parabolic systems by Hochbruck
and Ostermann \cite{HO}. In the stiff regime (here up to 
errors of order $10^{-4}$) DCRK
does not converge, the Integrating Factor method shows only first 
order behaviour (as predicted in \cite{HO}), the IRK4 scheme shows 
second order convergence,  and ETD methods perform 
best. This implies that the Cox-Matthews and Krogstad method with 
similar performance are the most economic for 
the stiff regime of the KP I equation, which gives the precision one 
is typically interested in in this context. For higher precisions we find
$a=1.87$ for the Integrating Factor method, $a=4.20$ for DCRK, $a=4.03$ for IRK4, 
$a=3.75$ for Krogstad's ETD scheme, $a=3.90$ for the Cox-Matthews scheme, and $a=3.96$ 
for the Hochbruck-Ostermann scheme.  

\begin{figure}[h!]
 \centering
 \includegraphics[width=0.49\textwidth]{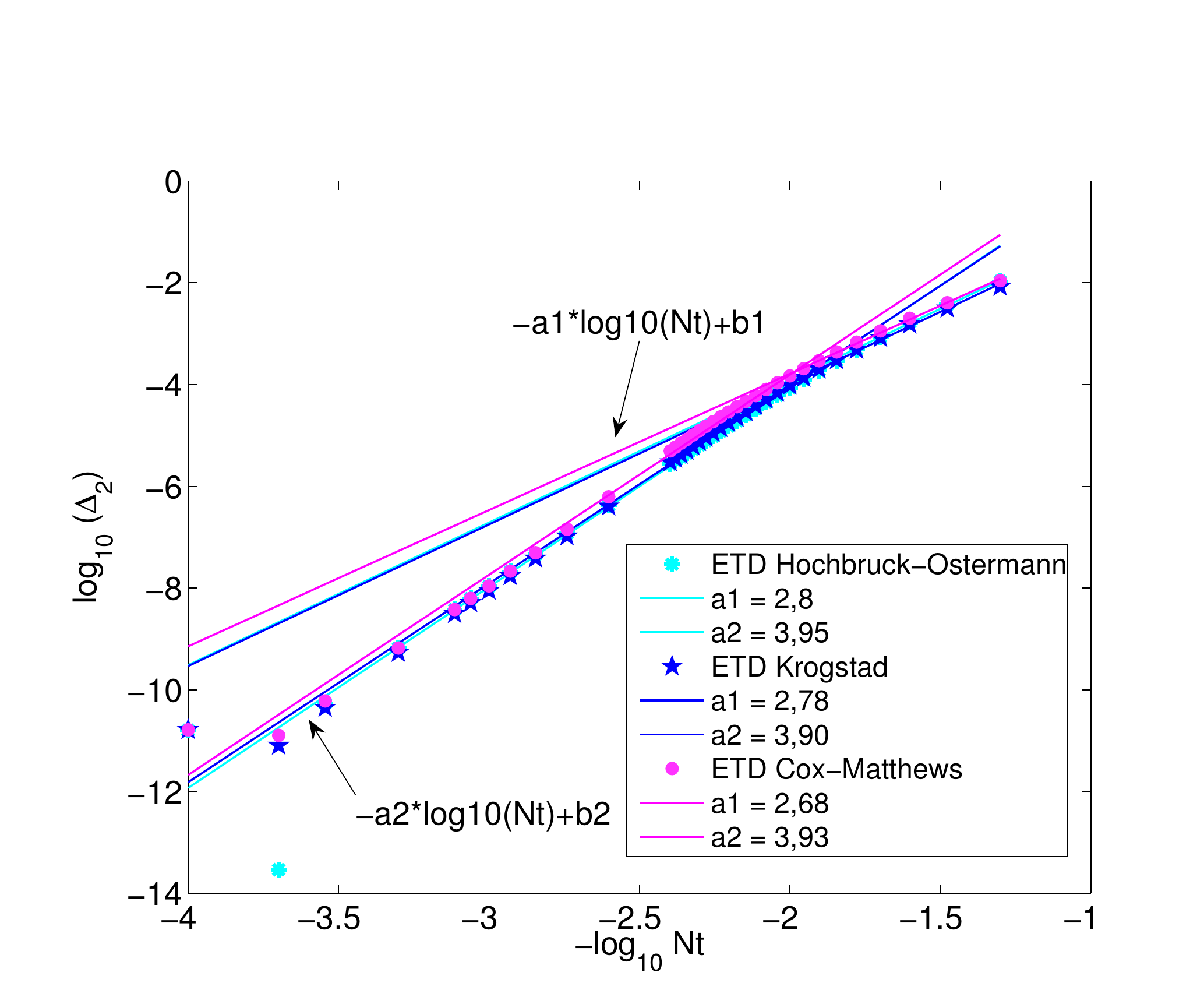}
 \includegraphics[width=0.49\textwidth]{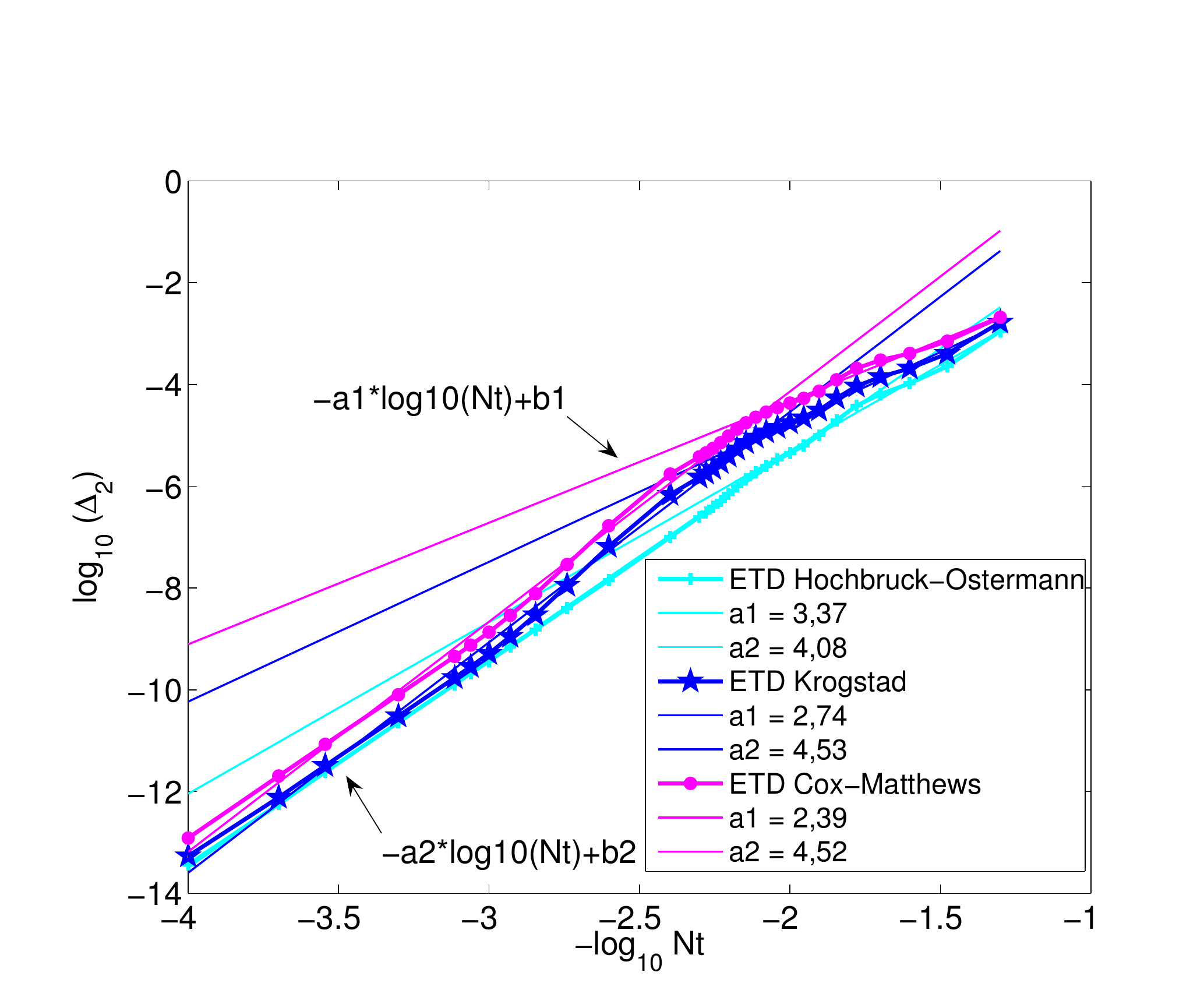}
\caption{Phenomenon of order reduction in exponential integrators 
for KP I  and KP II respectively in the small dispersion limit, ETD schemes of 
Fig.~\ref{figKPIsmallreg} (left) and of
Fig.~\ref{figKPIIsmallreg} (right).}
 \label{figKPIsmallreg2}
\end{figure}

To study empirically the phenomenon of order reduction
in exponential integrators, and to observe the transition from a 
stiff to a non stiff regime we study the ETD schemes in more detail 
in Fig.~\ref{figKPIsmallreg2}. This is
indicated by the fact that ETD schemes are only of order three in this stiff region instead of order four.
It appears that all schemes show a slight order reduction though this 
is not the case for the Hochbruck-Ostermann method in the parabolic 
case. 

\subsection{Numerical solution of Cauchy problems for the KP II equation}
\paragraph{Doubly periodic solution of KP II}

\begin{figure}[htb!]
 \centering
\includegraphics[width=0.49\textwidth]{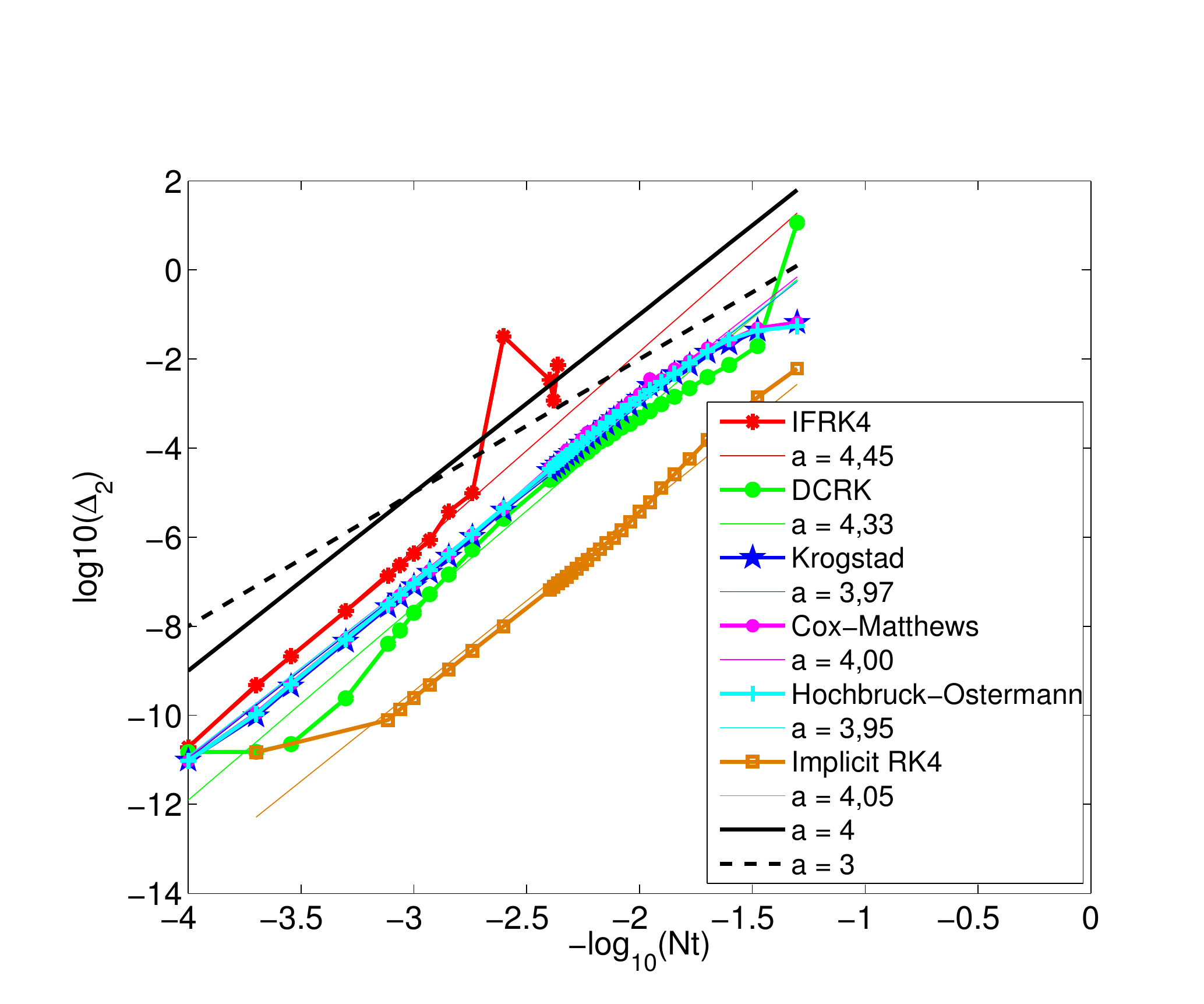}
\includegraphics[width=0.49\textwidth]{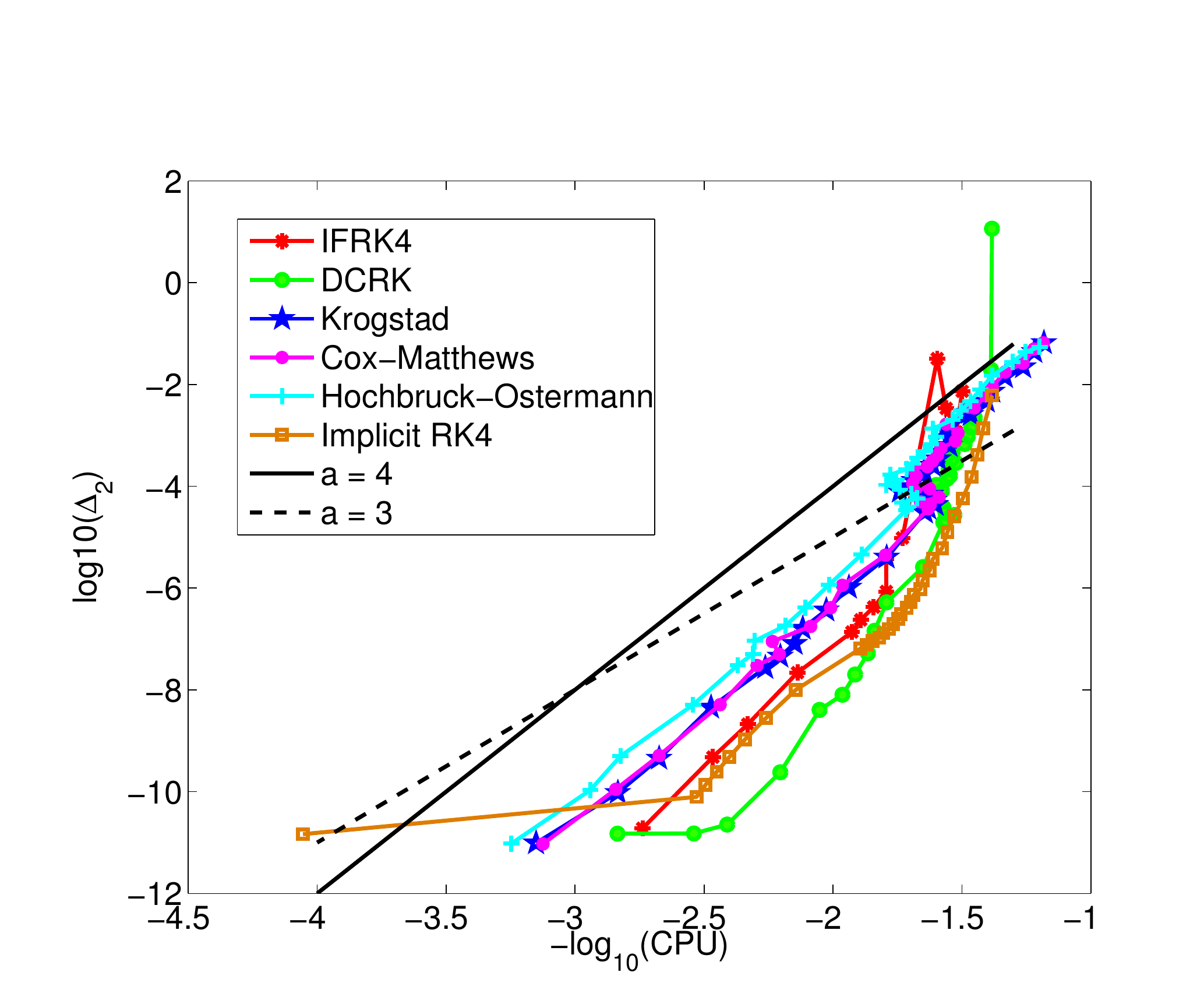}
  \caption{Normalized $L_{2}$ norm of the numerical error for the time 
 evolution of the doubly periodic solution to the KP II equation for several numerical
 methods, as a function of $N_{t}$  (left) and as a function of CPU 
 time (right).}
 \label{figperioreg}
\end{figure}

The computation for the doubly periodic solution to KP II 
is carried out with $2^{8}\times2^{8}$ points for $(x, y)\in[-5\pi,\,5\pi]\times[-5\pi,\,5\pi]$
and $t\leq1$ with the parameters $b = 1$, $\lambda = 0.15$, $b\lambda^{2}+d = -1$, $\mu_1 = \mu_2 = 0.25$, 
$ \nu_1 = -\nu_2 = 0.25269207053125$, $\omega_1 = \omega_2 = -1.5429032317052$,
and $\varphi_{1,0}=\varphi_{2,0} = 0$. 
The decrease of the numerical error is shown in 
Fig.~\ref{figperioreg} in dependence of $N_{t}$ and  in dependence on CPU time. 
From a linear regression analysis in a double logarithmic plot we can see that all schemes are
fourth order: one finds $a=4.45$ for the Integrating Factor method, $a=4.33$ for DCRK, 
$a=3.97$ for Krogstad's ETD scheme, $a=4$ for the Cox-Matthews 
scheme, and $a=3.95$ for the Hochbruck-Ostermann scheme.
As for the Zaitsev solution, DCRK performs best followed 
by the ETD schemes. We thus confirm Driscoll's results in \cite{D} 
on the efficiency of his method for this example. The absence of order reductions indicates again 
that the exact solution tests the equation in a non-stiff regime. 
IRK4 is competitive for larger time steps in this case since only 
very few iterations (1-3) are needed.

\paragraph{Small dispersion limit for KP II}

\begin{figure}[htb!]
 \centering
 \includegraphics[width=\textwidth]{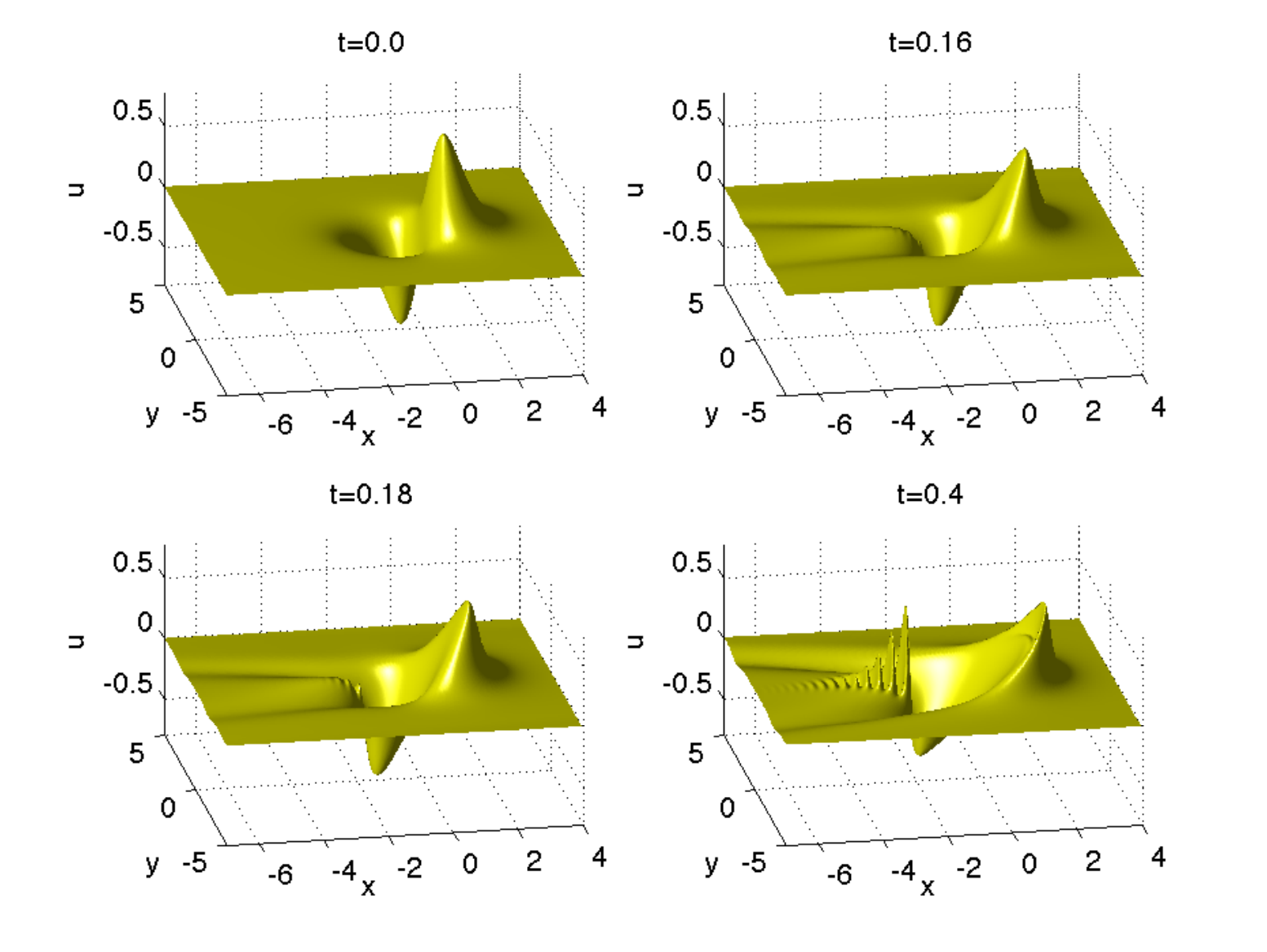}
 \caption{Solution to the KP II equation for the initial data $u_{0}=-\partial_{x}\mbox{sech}^{2}(R)\,\,\,\,\,\,
\mbox{where}\,\,\, R=\sqrt{x^{2}+y^{2}}$ for several values of $t$.}
    \label{figKPII}
\end{figure}

We consider the same initial data and the same methods as for KP I. 
In Fig.~\ref{figKPII} the time evolution of these data can be seen. 
The solution develops tails this time in negative $x$-direction. The 
steepening of the wave fronts happens at essentially the same time, 
but the gradients are stronger in the vicinity of the tails (see 
\cite{KSM}). This is also where the stronger oscillations appear. 

\begin{figure}[htb!]
   \centering
\includegraphics[width=0.49\textwidth]{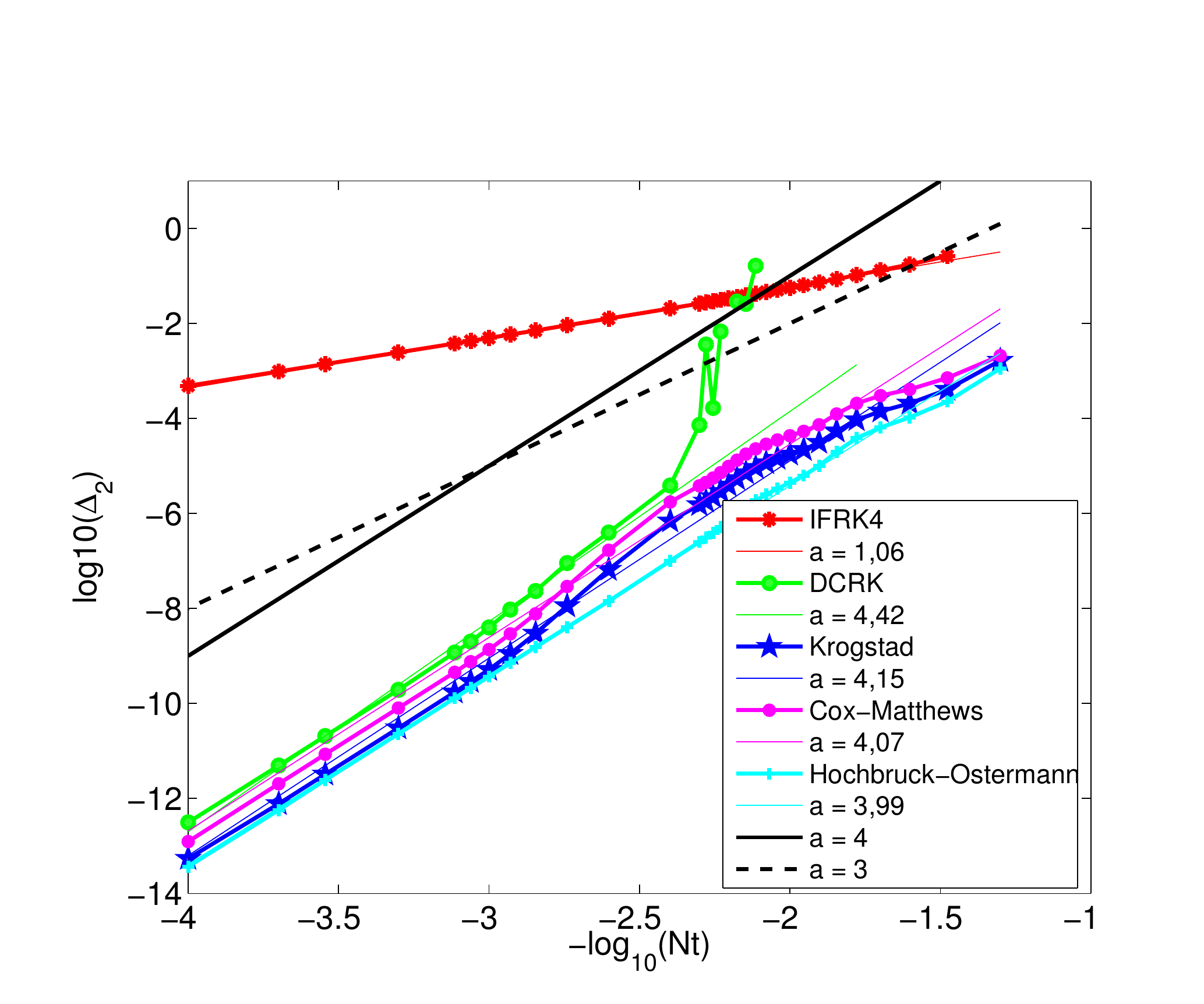}
\includegraphics[width=0.49\textwidth]{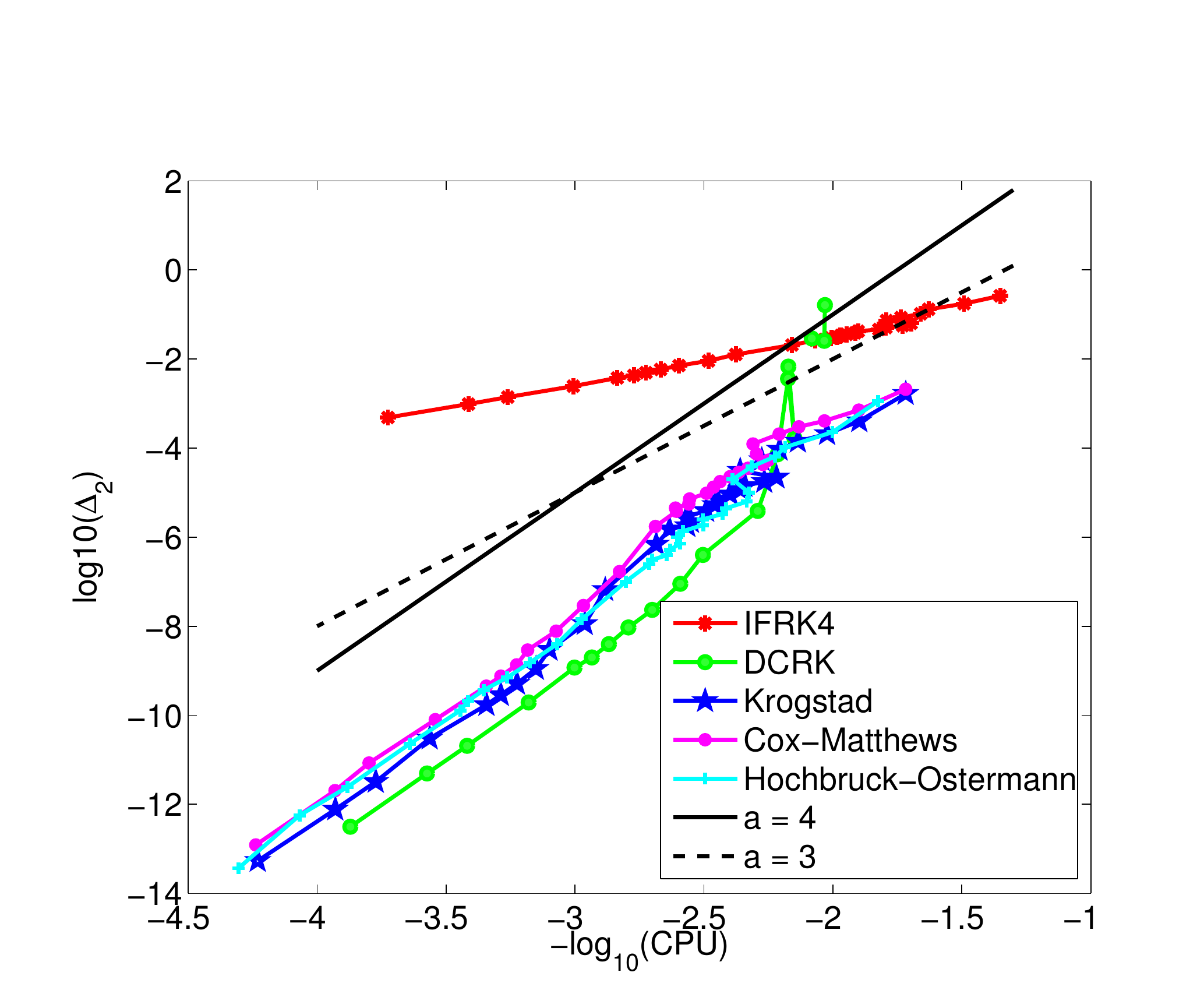}
 \caption{Normalized $L_{2}$ norm of the numerical error for the 
 solution in Fig.~\ref{figKPII} for several numerical
methods as a function of $N_{t}$ (left) and as a function of CPU 
 time (right).}
  \label{figKPIIsmallreg}
\end{figure}

The computation is carried out with $2^{11}\times2^{9}$
points for $(x, y)\in[-5\pi,\,5\pi]\times[-5\pi,\,5\pi]$,
$\epsilon=0.1$ and $t\leq0.4$. As a reference solution, we consider
the solution calculated with the Hochbruck-Ostermann method with $N_{t}=5000$
time steps. The dependence of the normalized $L_{2}$ norm of the
difference between this reference solution and the numerical solution
on $N_{t}$ and on CPU time is shown in Fig.~\ref{figKPIIsmallreg}.
 We obtain  similar results as in the 
small dispersion limit of KP I:
for typical accuracies one is interested in in this context, DCRK 
does not converge, the Integrating Factor method shows only a first 
order behavior, the IRK4 method is of second order, and the ETD methods perform best. In the non-stiff 
regime we find
 $a=1.06$ for the Integrating Factor method, $a=4.42$ for
 DCRK, $a=4.33$ for IRK4, 
$a=4.15$ for Krogstad's ETD scheme, $a=4.07$ for the Cox-Matthews scheme, and $a=3.99$ 
for the Hochbruck-Ostermann scheme.  
Once again, we study empirically the phenomenon of order reduction
in exponential integrators, and observe a transition from a stiff to a non 
stiff region (Fig.~\ref{figKPIsmallreg}), 
indicated by the fact that ETD schemes are only of order three in this stiff region instead of of order four.

\section{Davey-Stewartson equation}
In this section we perform a similar study as for KP of the efficiency of 
fourth methods in solving  Cauchy problems for the DS II equations. We first 
review some analytic facts about the focusing and defocusing DS II 
equations which are of importance 
in this context. 
We compare the performance of the codes for a 
typical example in the small dispersion limit.

\subsection{Analytic properties of the DS equations}
For a review see for instance the book by Sulem and 
Sulem \cite{SS}. 
We will study here only the DS II equations ($\alpha=\beta=1$ in 
eq.~(\ref{DSII})) since the 
elliptic operator for $\Phi$ can be inverted by imposing simple boundary 
conditions. For a hyperbolic operator acting on $\Phi$, boundary 
conditions for wave equations have to be used.

We will consider the equations again on  
$\mathbb{T}^{2}\times \mathbb{R}$. Due to the ellipticity of the 
operator in the equation for $\Phi$, it can be inverted in Fourier 
space in standard manner by imposing periodic boundary conditions on 
$\Phi$ as well. As before this case contains Schwartzian functions 
that are periodic for numerical purposes. Notice that solutions to 
the DS equations for Schwartzian initial data stay in this space at 
least for finite time in contrast to the KP case. Using Fourier 
transformations  $\Phi$ can be eliminated from the 
first equation by a transformation of the second equation in 
(\ref{DSII}) and an inverse transformation. With (\ref{fourier2d}) we have 
$$\Phi=-2\mathcal{F}^{-1}\left[\frac{k_{x}^{2}}{k_{x}^{2}+k_{y}^{2}}\mathcal{F}\left[|u|^{2}\right]\right],$$
which leads in (\ref{DSII}) as for KP to a nonlocal equation with a  Fourier 
multiplier. 
This implies that the DS 
equation requires an additional computational cost of 2 two-dimensional 
FFT per intermediate time step, thus doubling the cost with respect to the 
standard 2d NLS equation. Notice that from a numerical point of view 
the same applies to the elliptic-elliptic DS equation that is not 
integrable. Our experiments indicate that except for the additional 
FFT mentioned above, the numerical treatment of the 2d and higher 
dimensional NLS is analogous to the DS II  case studied here. The 
restriction to this case is entirely due to the fact that one can 
hope for an asymptotic description of the small dispersion limit in 
the integrable case. Thus we study initial data of the form 
$u_{0}(x,y)=a(x,y)\exp(ib(x,y)/\epsilon)$ with $a,b\in\mathbb{R}$, 
i.e., the semi-classical limit well known from the Schr\"odinger 
equation. Here we discuss only real initial data for convenience. 

It is known that DS solutions can have blowup. Results by Sung 
\cite{Sun} establish global existence in time for 
initial data $\psi_{0}\in L_{p}$, $1\leq p < 2$ with a Fourier 
transform $\mathcal{F}[\psi_{0}]\in L_{1}\cap L_{\infty}$ subject to the 
smallness condition 
\begin{equation}
    ||\mathcal{F}[\psi_{0}]||_{L_{1}}||\mathcal{F}[\psi_{0}]||_{L_{\infty}}
    <\frac{\pi^{3}}{2}\left(\frac{\sqrt{5}-1}{2}\right)^{2}
    \label{sungcond}
\end{equation}
in the 
focusing case. There is no such condition in the defocusing case. 
Notice that condition (\ref{sungcond}) has been established for the DS II 
equation with $\epsilon=1$. The coordinate change $x'=x/\epsilon$, 
$t'=t/\epsilon$ transforms the DS equation (\ref{DSII}) to this 
standard form. 
This implies for the initial data $u_{0}=\exp(-x^{2}-\eta y^{2})$ we study for the small dispersion 
limit of the focusing DS II system  in this paper that condition 
(\ref{sungcond}) takes the form 
$$\frac{1}{\epsilon^{2}\eta}\leq 
\frac{1}{8}\left(\frac{\sqrt{5}-1}{2}\right)^{2} \sim 0.0477.$$
This condition is not satisfied for the values of $\epsilon$ and 
$\eta$ we use here. Nonetheless we do not observe any indication 
of blowup on the shown timescales. One of the reasons is that the 
rescaling with $\epsilon$ above also rescales the critical time for 
blowup by a factor $1/\epsilon$. In addition it is expected 
that the dispersionless equations will for generic initial data have 
a gradient catastrophe at some time $t_{c}<\infty$, and that the 
dispersion will regularize the solution for small times $t>t_{c}>0$. 
However there are no analytic results in this context.

The complete integrability of the DS II equation implies again the 
existence of explicit solutions. Multi-soliton solutions will be as in 
the KP case localized in one spatial direction and infinitely 
extended in another, the lump solution is localized 
in two spatial directions, but with an algebraic fall off towards 
infinity. Thus these are again not convenient to test codes 
based on Fourier methods as in the KP case.  Since the study of 
the small dispersion limit below indicates that the time steps 
have to be chosen sufficiently small for accuracy reasons such that 
in contrast to KP no order reduction observed, we will not 
study any exact solutions here.

\subsection{Small dispersion limit for DS II in the defocusing case}

\begin{figure}[htb!]
 \centering
 \includegraphics[width=\textwidth]{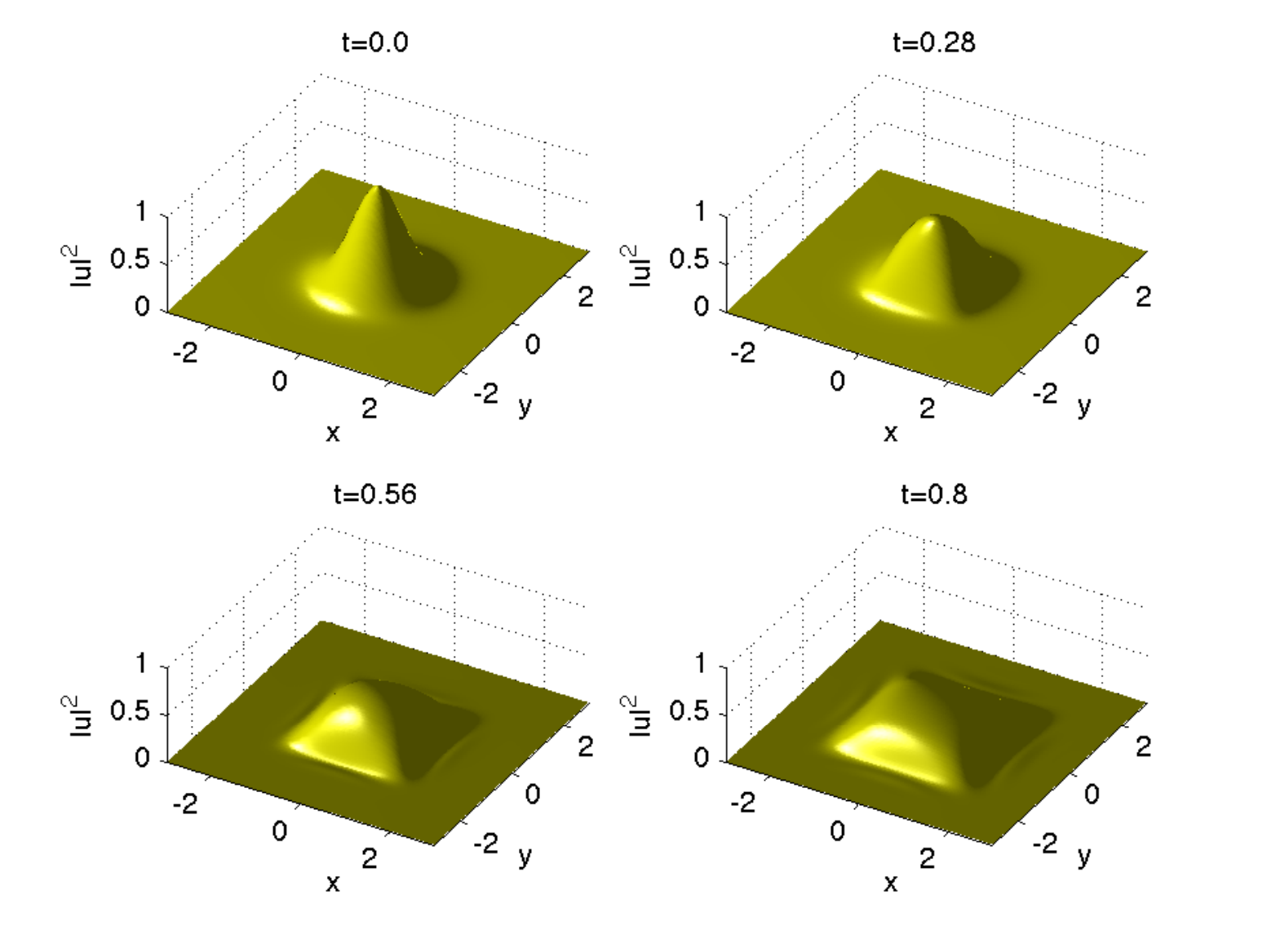}
 \caption{Solution to the defocusing DS II equation for the initial 
 data $u_{0}=\exp(-R^{2})$
 where  $R=\sqrt{x^{2}+y^{2}}$  and $\epsilon=0.1$ for several values of $t$.}
    \label{figDSdefoc}
\end{figure}

We consider initial data $u_{0}$ of the form
\begin{equation}\label{e13}
u_{0}(x,\, y)=e^{-R^{2}},\,\,\,\mbox{where}\,\,\, R=\sqrt{x^{2}+\eta y^{2}}\,\,\,\mbox{with}\,\,\,\eta=1
\end{equation}
and use the same methods as before together with time splitting 
methods of order 2 and one of order 4, as 
explained in section 2.4. The defocusing effect of the defocusing DS 
II equation for these initial data can be seen in 
Fig.~\ref{figDSdefoc}, where $|u|^{2}$ is shown for several values 
of $t$. The compression of the initial pulse 
into some almost pyramidal shape
leads to a steepening on the 4 sides parallel to the coordinate axes 
and to oscillations in these regions.

\begin{figure}[htb!]
 \centering
 \includegraphics[width=0.49\textwidth]{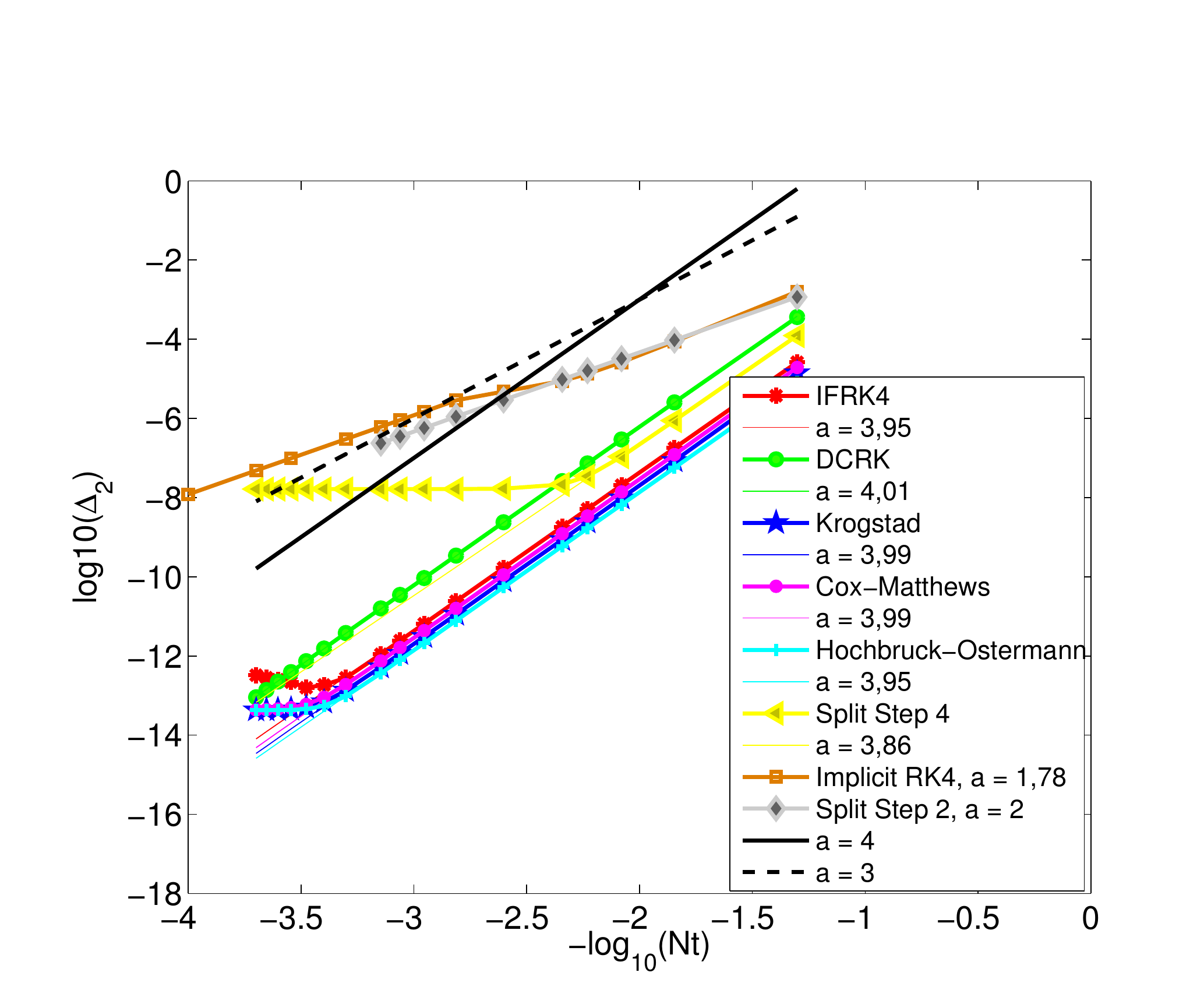}
 \includegraphics[width=0.49\textwidth]{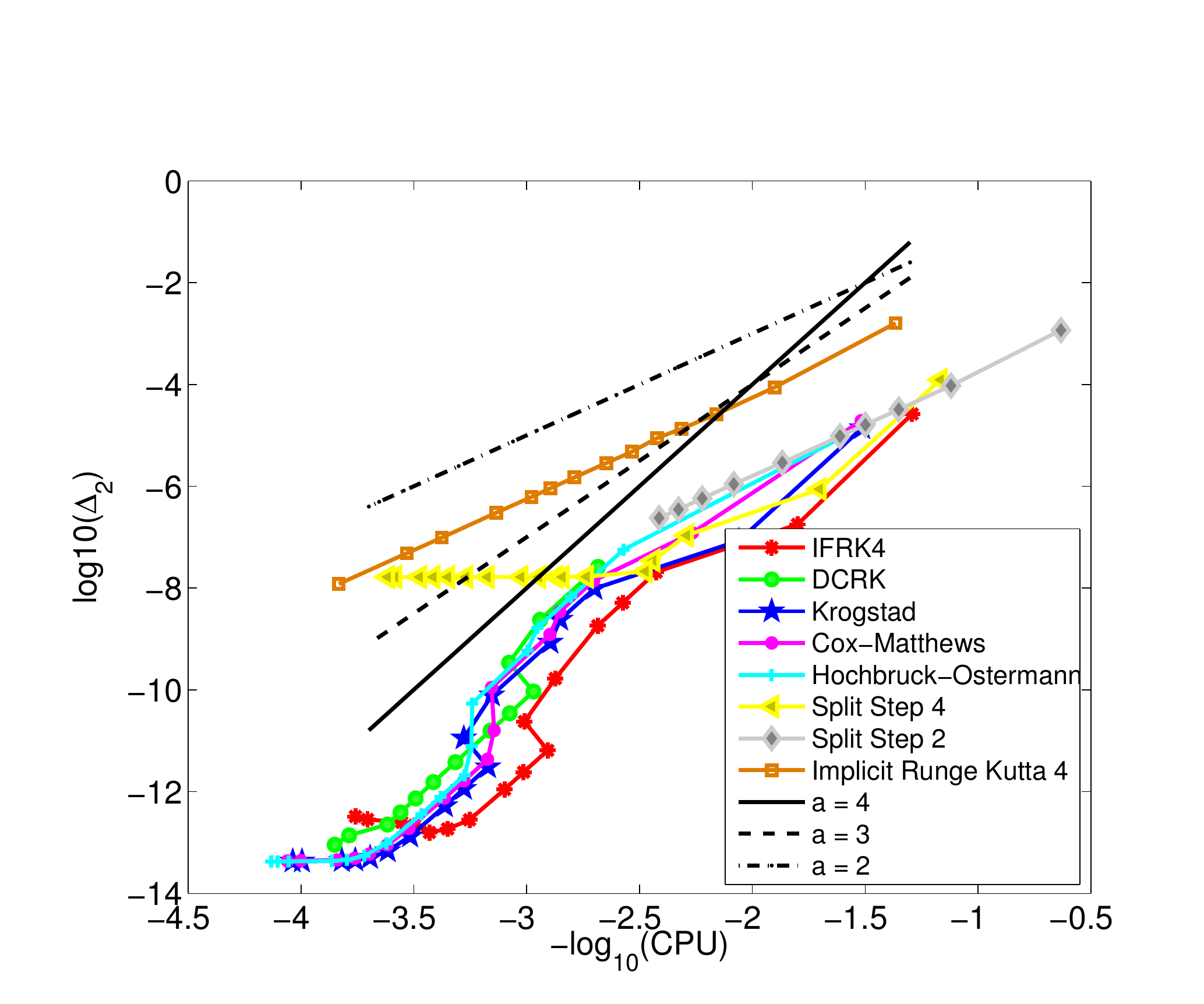}
 \caption{Normalized $L_{2}$ norm of the numerical error for several numerical
methods for the situation shown in Fig.~\ref{figDSdefoc} as a 
function of $N_{t}$ (left) and of CPU time (right).}
    \label{figDSdefocsmallreg}
\end{figure}

The computations are carried out with $2^{10}\times2^{10}$ points
for $(x, y)\in[-5\pi,\,5\pi]\times[-5\pi,\,5\pi]$, $\epsilon=0.1$
and $t\leq0.8$. To determine a reference solution, we compute solutions
with $6000$ time steps with the ETD, the DCRK and the IF schemes
and take the arithmetic mean. The dependence of the normalized $L_{2}$
norm of the difference of the numerical solutions with respect to
this reference solution on $N_{t}$
and on CPU time is shown in 
Fig.~\ref{figDSdefocsmallreg}.
A linear regression shows
 that all fourth order schemes show a fourth order behavior except 
for IRK4 ($a=1.78$), as is obvious from the 
 straight lines with slope
 $a=3.95$ for the Integrating Factor method, $a=4.01$ for
 DCRK,
$a=3.99$ for Krogstad's ETD scheme, $a=3.99$ for the Cox-Matthews 
scheme, $a=3.95$ 
for the Hochbruck-Ostermann scheme, and $a=3.86$ for the time 
splitting method. The second order splitting scheme shows the 
expected convergence rate and performs very well for lower precision. For 
smaller time steps, the advantage of the fourth order schemes is more 
pronounced.  Apparently the system is `stiff' for the IRK4 
scheme since it only shows second order behavior.
Notice that the time splitting scheme reaches its maximal precision 
around $10^{-8}$, a behavior which was 
already noticed in \cite{ckkdvnls} for the study 
of the Nonlinear Schr\"odinger equation
 in the small dispersion limit. It appears that this behavior is due 
 to resonances of errors of the split equations, but the 
 identification of the precise 
 reason will be the subject of further research. The same effect is observed for second 
 order splitting for smaller time steps than shown in 
 Fig.~\ref{figDSdefocsmallreg}. However, both schemes work very well 
 at the precisions in which one is normally interested in. 
 We mainly include a second order 
 scheme here because of the additional computational cost due to the 
 function $\Phi$ in the DS system. This could make a second order 
 scheme competitive in terms of CPU time because of the lower number 
 of FFT used per time step. It can be seen in 
 Fig.~\ref{figDSdefocsmallreg} that this is not the case.
We conclude that the ETD schemes perform best in this context.

\subsection{Small dispersion limit for the focusing DS II equation}

\begin{figure}[htb!]
  \centering
 \includegraphics[width=\textwidth]{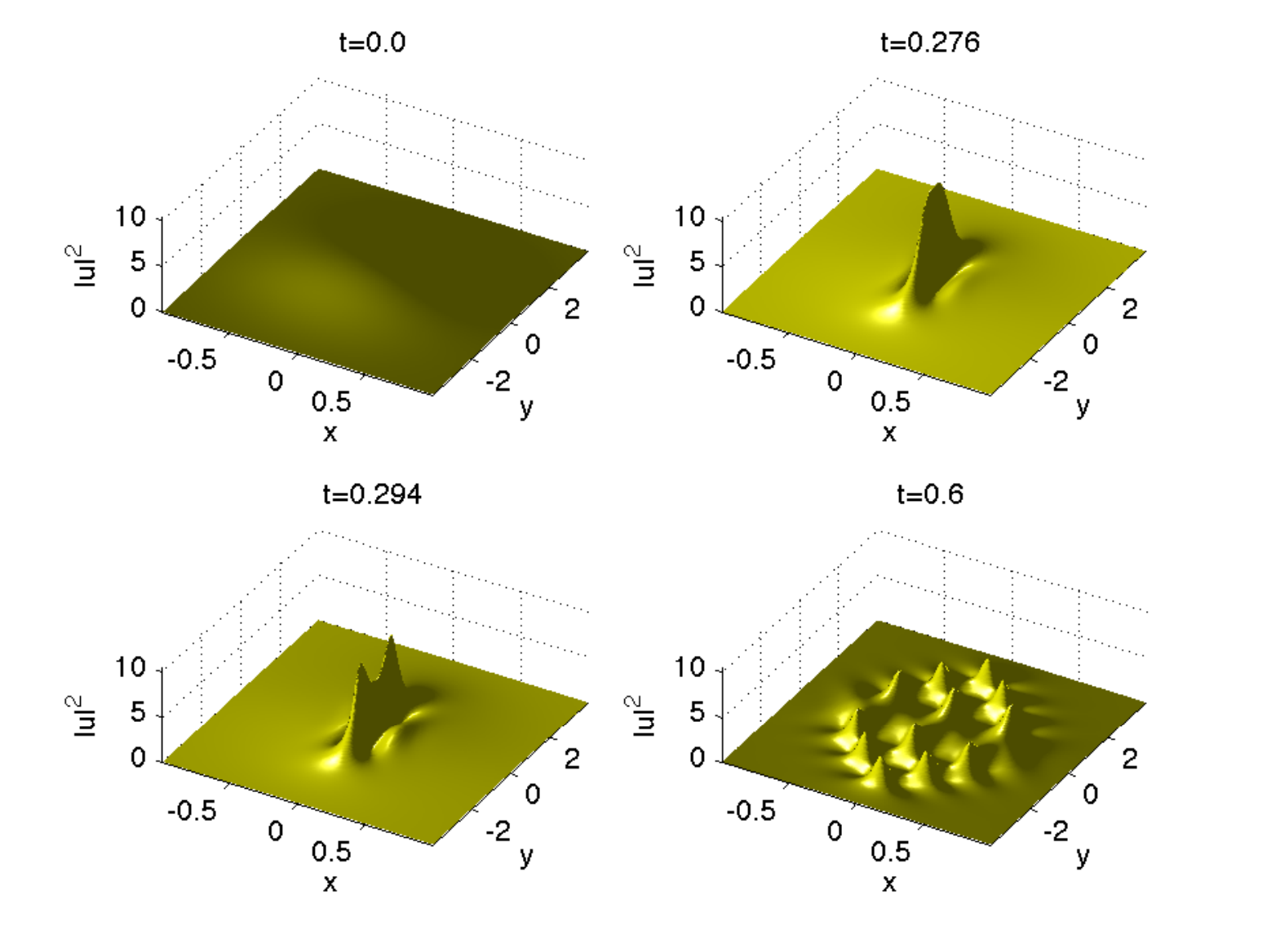}
 \caption{Solution to the focusing DS II equation for the initial 
 data $u_{0}=\exp(-R^{2})\,\,\,\,\,\,
\mbox{where}\,\,\, R=\sqrt{x^{2}+\eta y^{2}}$, $\eta=0.1$  and $\epsilon=0.1$ for several values of $t$.}
   \label{figDSfoc}
\end{figure}

For the focusing DS II in the small dispersion limit we consider initial data of the form
(\ref{e13}) with $\eta=0.1$
and the same methods as before. The focusing effect of the equation 
can be clearly recognized  in Fig.~\ref{figDSfoc}. The initial peak 
grows until a breakup into a pattern of smaller peaks occurs.

\begin{figure}[htb!]
  \centering
 \includegraphics[width=0.49\textwidth]{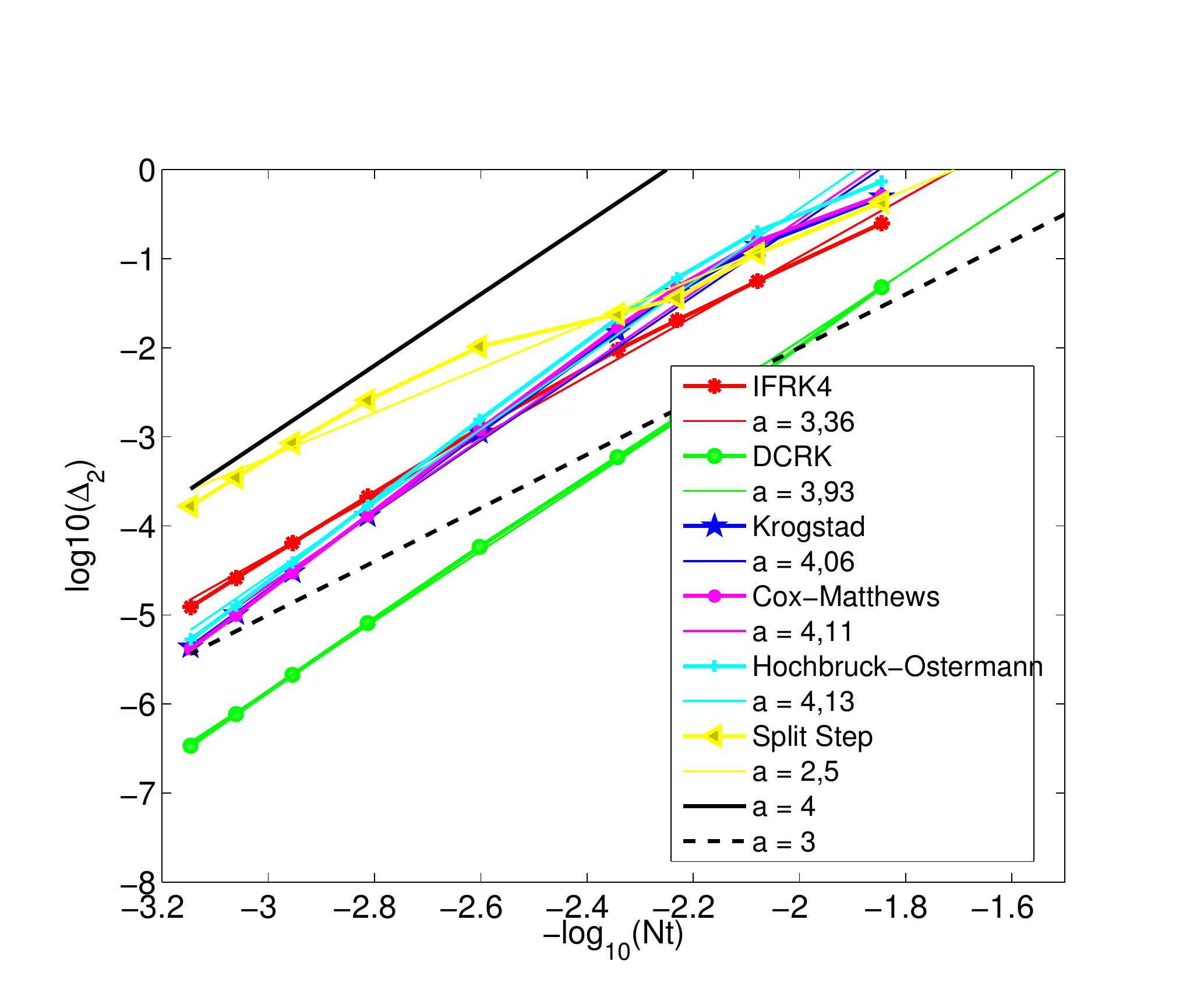}
 \includegraphics[width=0.49\textwidth]{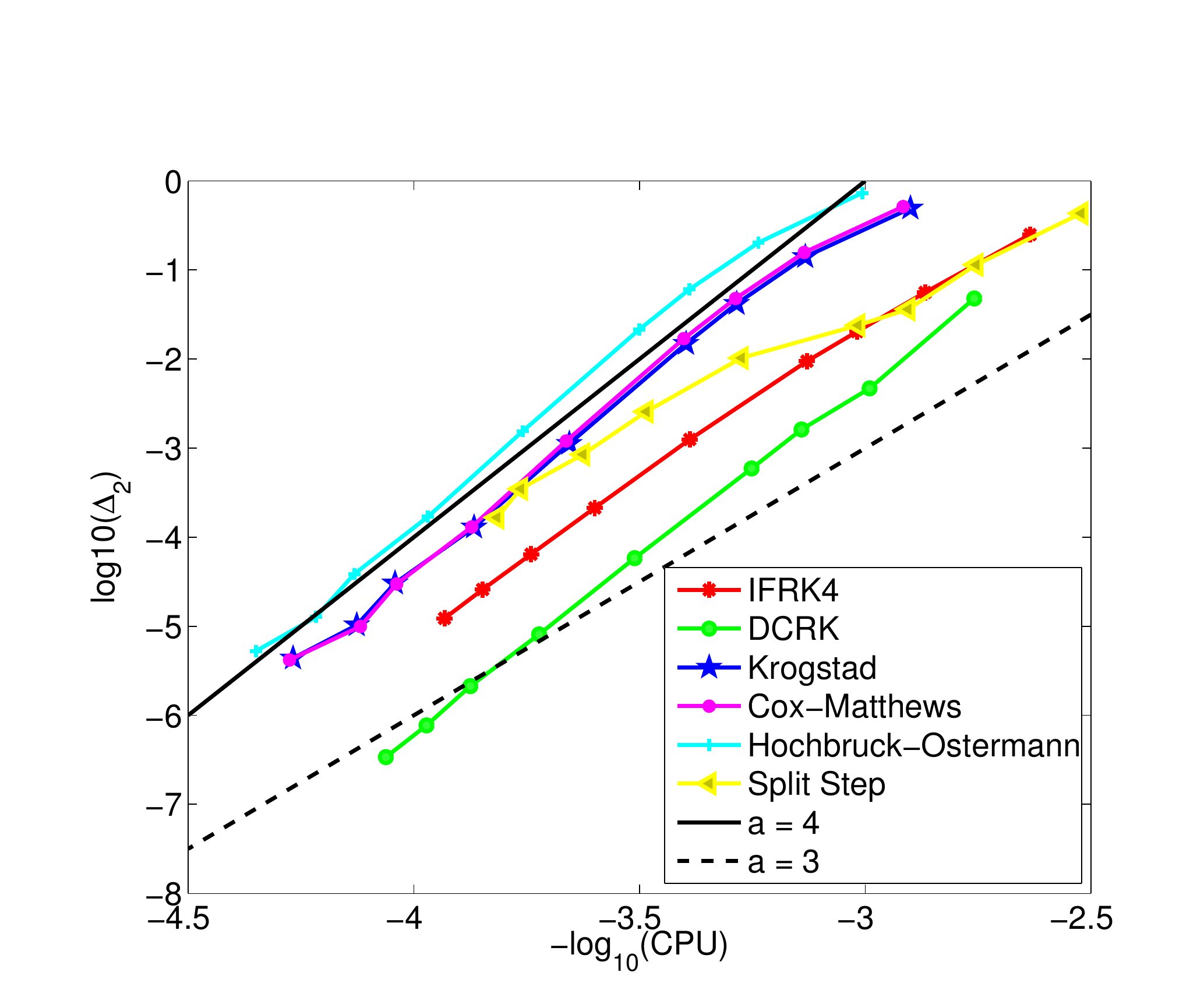}
 \caption{Normalized $L_{2}$ norm of the numerical error for the 
 example of Fig.~\ref{figDSfoc} for several numerical
methods as a function of $N_{t}$ (left) and of CPU time (right).}
   \label{figDSfocsmallreg}
\end{figure}

It is crucial to provide sufficient 
spatial resolution for the central peak. As for the 
1+1-dimensional focusing NLS discussed in \cite{ckkdvnls}, the 
modulational instability of the focusing DS II leads to numerical 
problems if there is no sufficient resolution for the maximum.
In \cite{ckkdvnls} a 
resolution of $2^{13}$ modes was necessary for initial data 
$e^{-x^{2}}$ and $\epsilon=0.1$ for the focusing NLS in 
$1+1$ dimensions.  The possibility of blowup in DS requires at 
least the same resolution despite some regularizing effect of the 
nonlocality $\Phi$. With the computers we could access, a systematic 
study of time integration schemes with a resolution of $2^{13}\times 
2^{13}$ was not possible in Matlab. Thus we settled for initial data close to 
the one-dimensional case, which allowed for a lower resolution, see 
the next section for the Fourier coefficients.
The computation is carried out with $2^{12}\times2^{11}$ points
for $(x, y)\in[-5\pi,\,5\pi]\times[-5\pi,\,5\pi]$, $\epsilon=0.1$
and $t\leq0.6$. To determine a reference solution, we compute solutions
with $6000$ time steps with the ETD, the DCRK and the IF schemes
and take the arithmetic mean. The dependence of the normalized $L_{2}$
norm of the difference of the numerical solutions with respect to
this reference solution on 
$N_{t}$ and on CPU time is shown in Fig.~\ref{figDSfocsmallreg}.
All schemes except the time splitting scheme show a fourth order 
 behavior, as can be seen from the straight lines with slope
 $a=3.36$ for the Integrating Factor method, $a=3.93$ for
 DCRK, 
$a=4.06$ for Krogstad's ETD scheme, $a=4.11$ for the Cox-Matthews scheme, $a=4.13$ 
for the Hochbruck-Ostermann scheme, and $a=2.5$ for the fourth order
time splitting method. 
We conclude that in this context  DCRK performs 
best, followed by the ETD schemes. We do not present results for the 
IRK4 scheme here since it was computationally too expensive.

\section{Numerical conservation of the $L_{2}$ norm}
The complete integrability of the KP and the DS equations implies the 
existence of many or infinitely many conserved quantities (depending 
on the function spaces for which the solutions are defined). It can 
be easily checked that the $L_{1}$ norm and the $L_{2}$ norm of the 
solution are conserved as well as the energy. We do not use here 
symplectic integrators that take advantage of the Hamiltonian 
structure of the equations. Such integrators of fourth order will be 
always implicit which will be in general computationally too expensive for the 
studied equations as the experiment with the implicit IRK4 scheme showed. 
Moreover it was shown in \cite{BIS} that fourth order exponential 
integrators clearly outperform second order symplectic integrators 
for the NLS equation. 

The fact that the conservation of $L_{2}$ norm and energy is not 
implemented in the code allows to use the `numerical conservation' of 
these quantities during the computation or 
the lack thereof to test the quality of the 
code. We will study in this section for the previous examples to 
which extent this leads to a quantitative 
indicator of numerical errors. Note that due to the non-locality of 
the studied PDEs (\ref{e1}) and (\ref{DSII}), the energies both for KP,
$$
   E[u(t)]:= \frac{1}{2} \int_{\mathbb{T}^2} \left(\partial_x u(t,x,y))^2 - \lambda 
   (\partial_x^{-1} \partial_y u(t,x,y))^2 - 2\epsilon^2 u^3(t,x,y) 
   \right)
   d x d y,$$
and for DS II,
\begin{align}
    E[u(t)] & := \frac{1}{2} \int_{\mathbb{T}^2} \bigg[ \epsilon^{2}|\partial_x u(t,x,y)|^2 - 
\epsilon^{2}|\partial_y u(t,x,y)|^2     \nonumber\\
& \left.-\rho\left(|u(t,x,y)|^{4}-\frac{1}{2}\left(\Phi(t,x,y)^{2}+(\partial_{x}^{-1}\partial_{y}\Phi(t,x,y))^{2}\right)\right) 
  \right] d x d y,
    \nonumber
\end{align}
contain anti-derivatives with respect to $x$. Since the latter are 
computed with Fourier methods, i.e., via division by  $k_{x}$ in 
Fourier space, this 
computation is in itself numerically problematic and could indicate 
problems not present in the numerical solution of the Cauchy problem. Therefore we trace 
here only the $L_{2}$ norm $\int_{\mathbb{T}^2}|u(t,x,y)|^{2}dxdy $, 
where these problems do not appear.  In the plots we show the 
variable $test$ defined as $test=M(t)/M(0)-1$, where $M(t)$ is the 
numerically computed $L_{2}$ norm in dependence of time.

Notice that numerical conservation of the $L_{2}$ norm can be only 
taken as an indication of the quality of the numerics if there is 
sufficient spatial resolution. Therefore we will always present the 
Fourier coefficients for the final time step for the considered 
examples. No dealiasing techniques are used.
We will discuss below the results for the small dispersion 
limit. 

\begin{figure}[htb!]
    \includegraphics[width=0.5\textwidth]{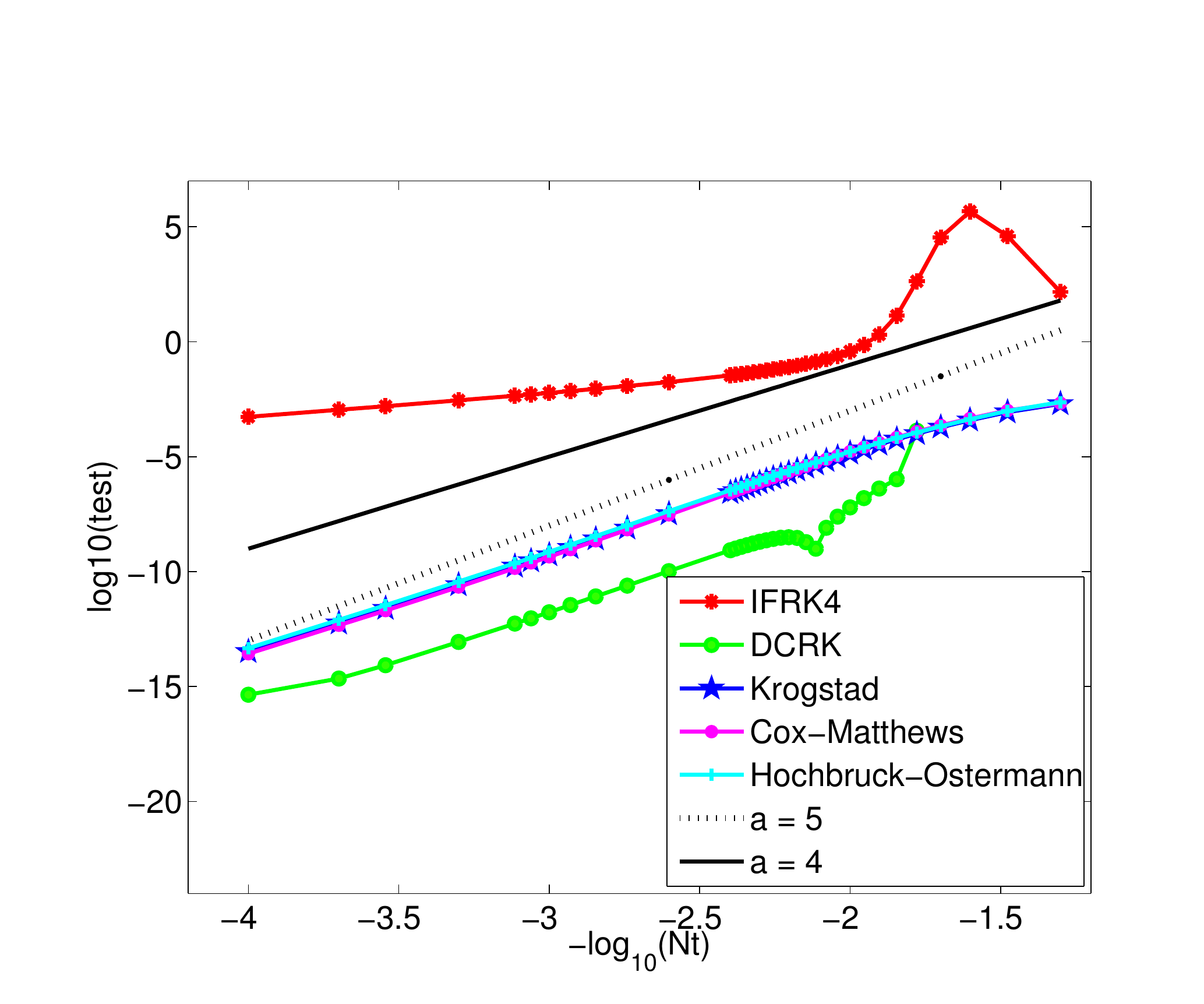}
 \includegraphics[width=0.5\textwidth]{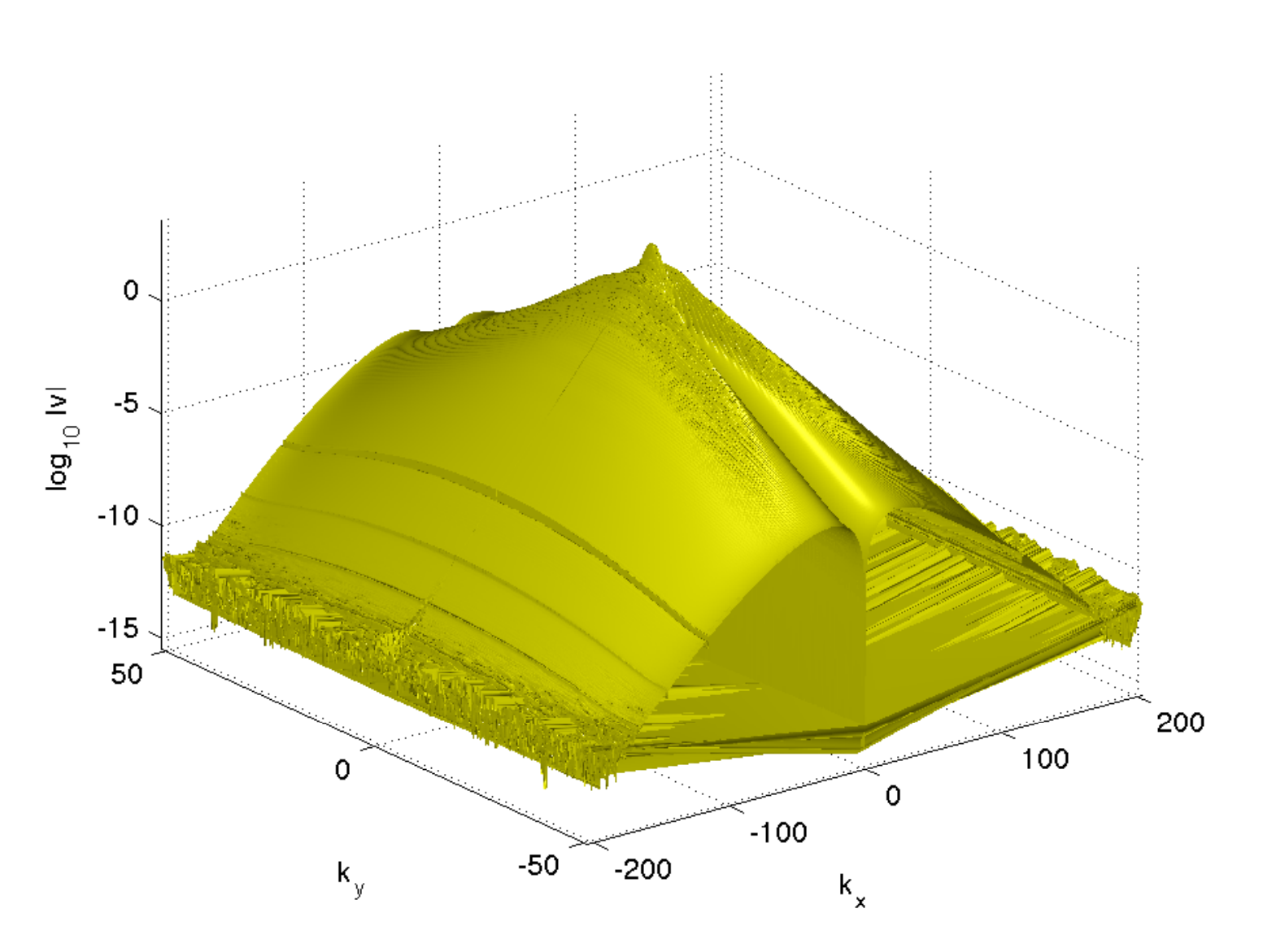}
 \caption{Non-conservation of the numerically computed $L_{2}$ norm 
 of the solution to the problem considered in Fig.~\ref{figKPIsmallreg} 
 in dependence on the time step
(left) and the Fourier coefficients for the final time 
 (right). }
   \label{figKPIsmallmass}
\end{figure}

For the KP I example of Fig.~\ref{figKPI} we get the Fourier 
coefficients at the final time and the mass conservation shown in 
Fig.~\ref{figKPIsmallmass}. It can be seen that the Fourier 
coefficients decrease in $k_{x}$-direction to almost machine 
precision, whereas this is not fully achieved in $k_{y}$-direction. 
This is partly due to the necessity to allow extensive studies of the 
dependence on the time-stepping in finite computing time and thus to 
keep the spatial resolution low, and partly due 
to a Gibbs phenomenon mainly in $k_{y}$-direction due to the formation of the 
algebraic tails in Fig.~\ref{figKPI}. Mass conservation can be seen to 
be a viable indicator of the numerical accuracy by comparing with 
Fig.~\ref{figKPIsmallreg}: in the range of accuracy in which one is 
typically interested ($\sim10^{-4}$), mass conservation overestimates the actual 
accuracy by roughly 2 orders of magnitude. It can be seen that it 
shows also at least a fourth order decrease. 

\begin{figure}[htb!]
 \includegraphics[width=0.5\textwidth]{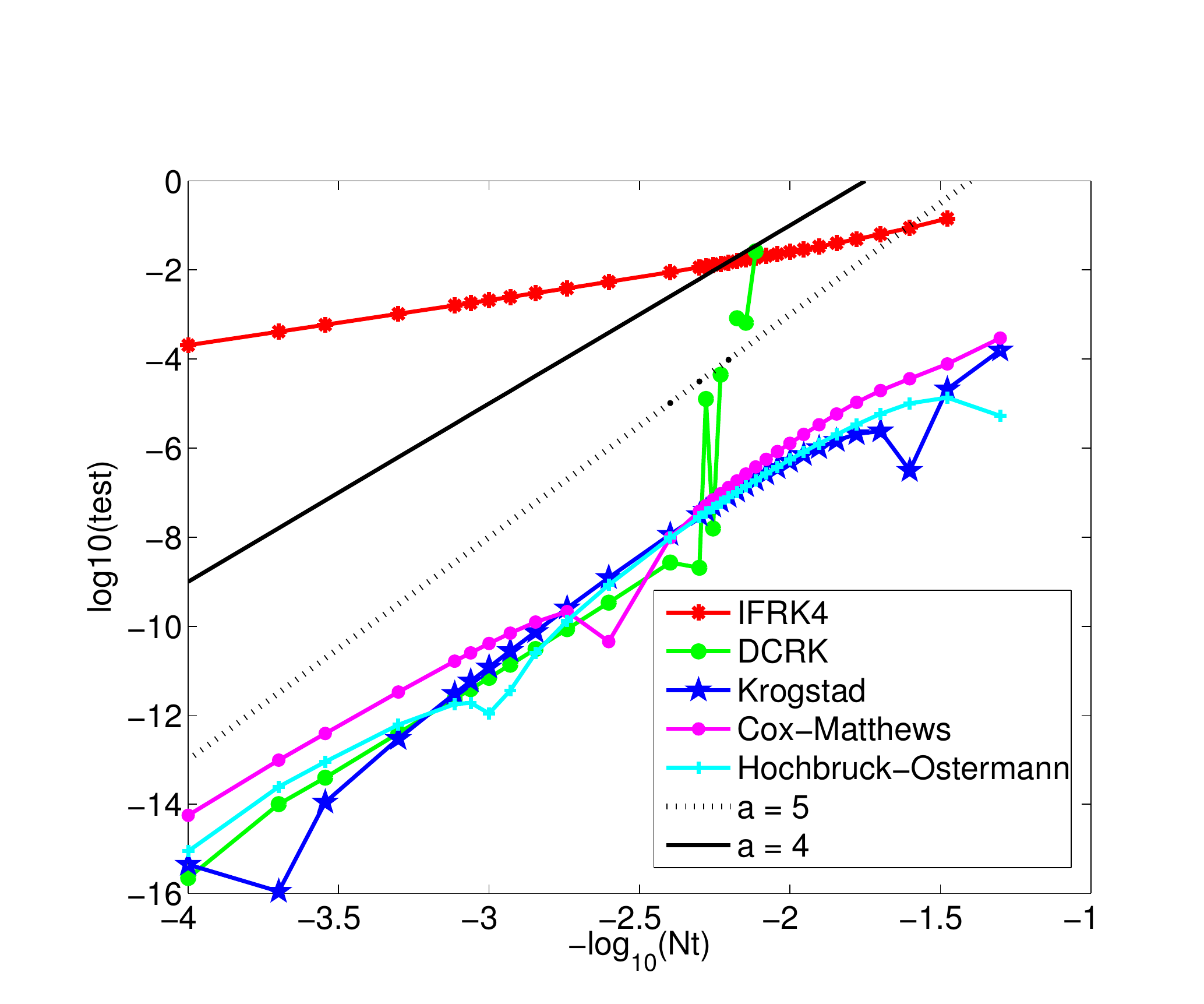}
 \includegraphics[width=0.5\textwidth]{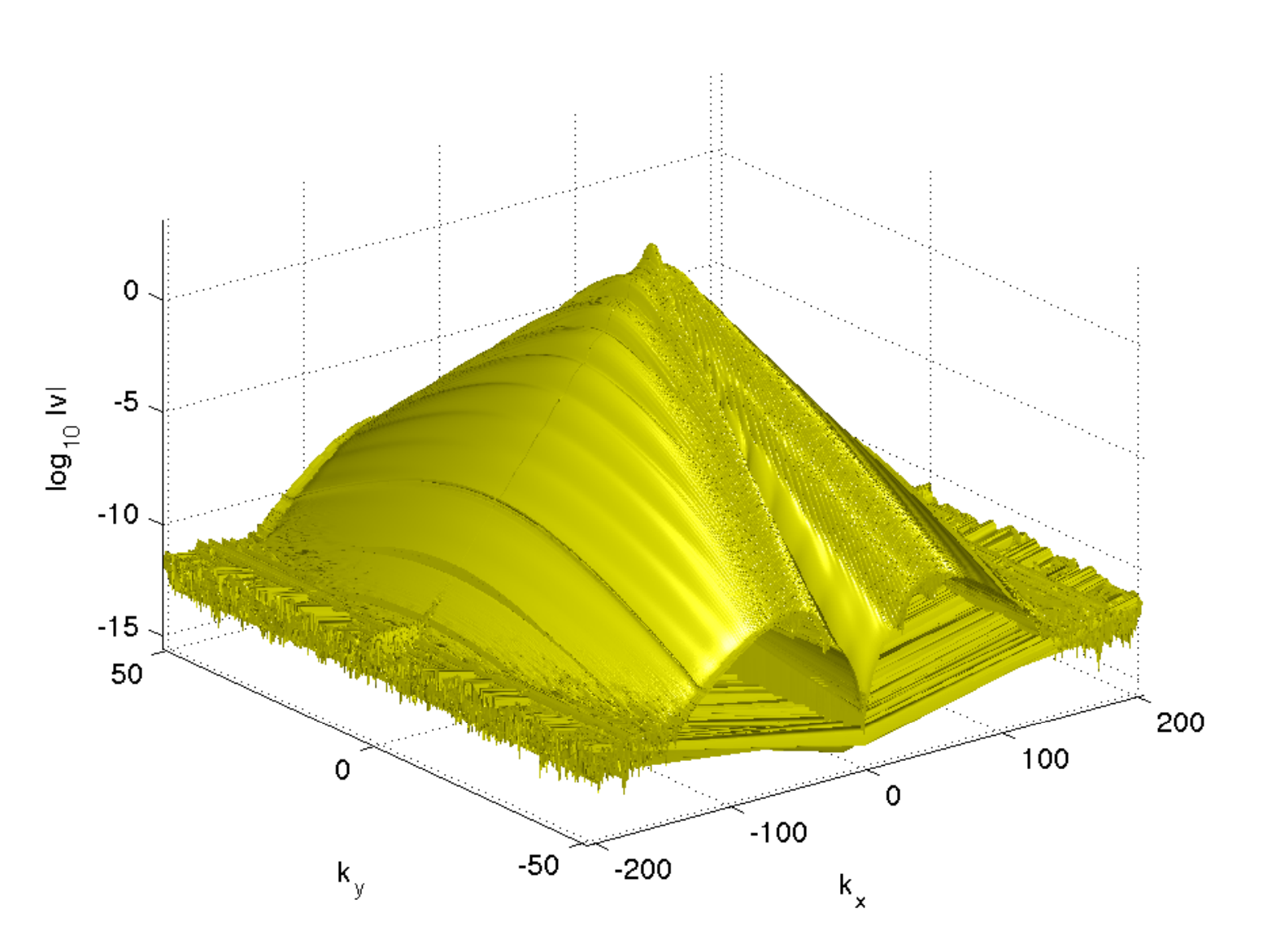}
 \caption{Non-conservation of the numerically computed $L_{2}$ norm 
 of the solution to the problem considered in 
 Fig.~\ref{figKPIIsmallreg} 
 in dependence on the time step
 (left) and the Fourier coefficients for the final time 
  (right). }
    \label{figKPIIsmallmass}
\end{figure}

The situation is very similar for the small dispersion example for KP II 
of Fig.~\ref{figKPII} as can be seen in Fig.~\ref{figKPIIsmallmass}.

For the defocusing DS II equation and the example shown in 
Fig.~\ref{figDSdefoc}, the Fourier coefficients decrease to machine 
precision despite the lower resolution than for KP. One reason 
for this is the absence of algebraic tails in the solution. The mass shows as for KP 
at least fourth order dependence on the time step and overestimates the 
numerical precision by roughly two orders of magnitude. This is not 
true for the splitting scheme for which mass conservation is no 
indication of the numerical precision at all. This seems to be due to 
the exact integration of the equations (\ref{2.4}) into which DS is 
split (for one of them the $L_{2}$ 
norm is constant). The found numerical mass does not appear to reflect the 
splitting error that is the reason for the numerical error here. 
\begin{figure}[htb!]
 \includegraphics[width=0.5\textwidth]{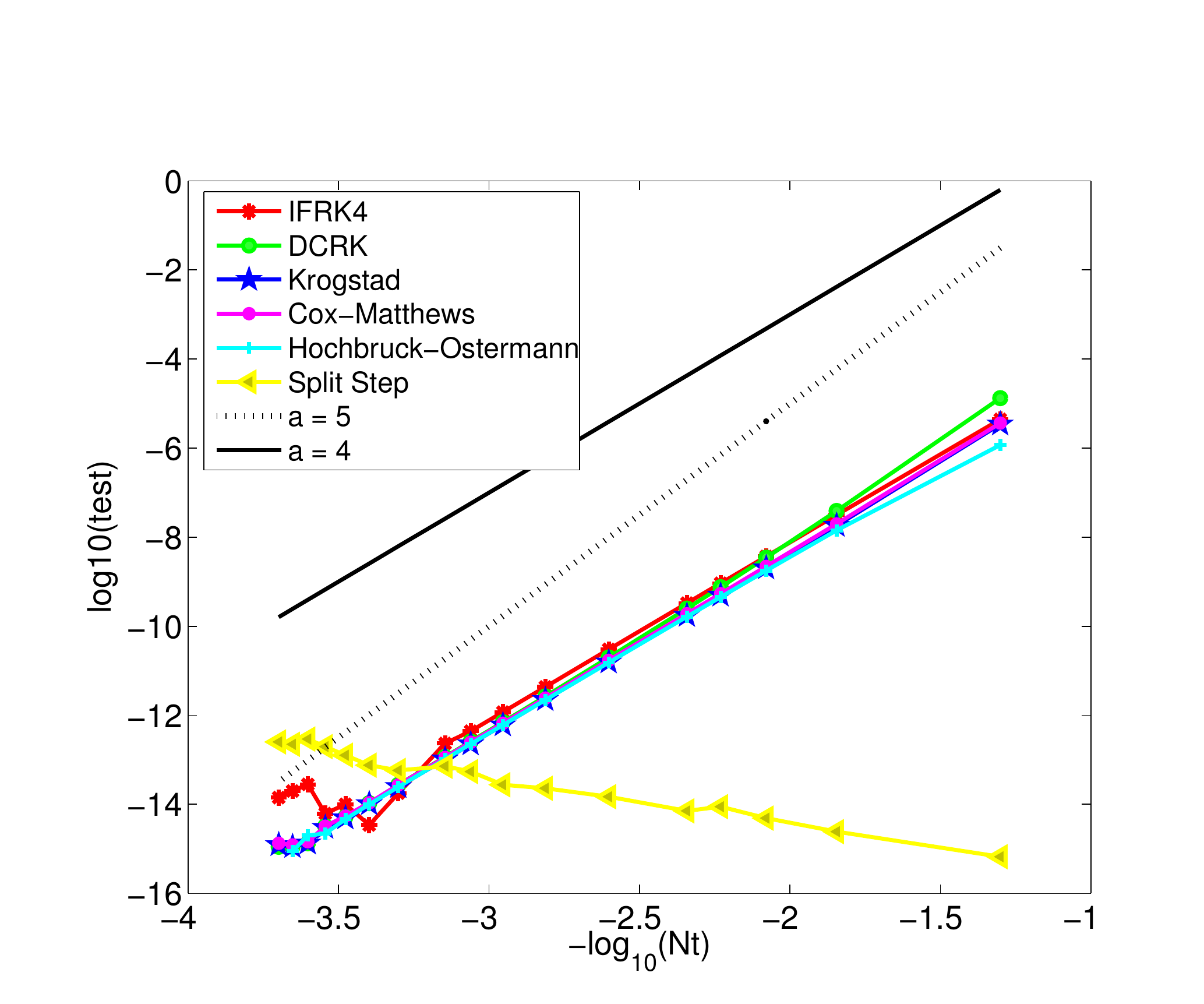}
 \includegraphics[width=0.5\textwidth]{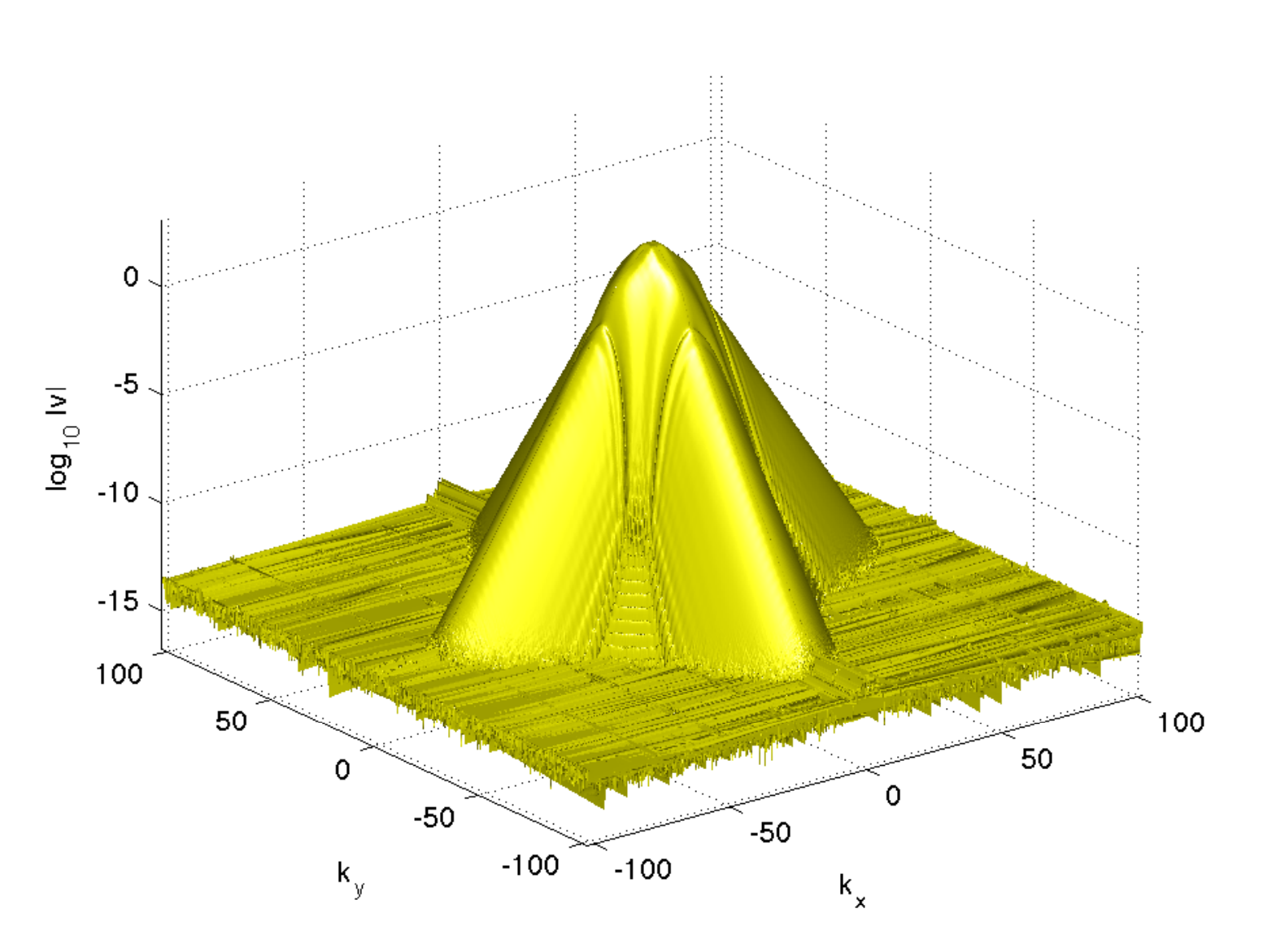}
 \caption{Non-conservation of the numerically computed $L_{2}$ norm 
 of the solution to the problem considered in 
 Fig.~\ref{figDSdefocsmallreg} 
 in dependence on the time step
 (left) and the Fourier coefficients for the final time 
  (right). }
   \label{figDSdefocsmallmass}
\end{figure}

\begin{figure}[htb!]
  \includegraphics[width=0.5\textwidth]{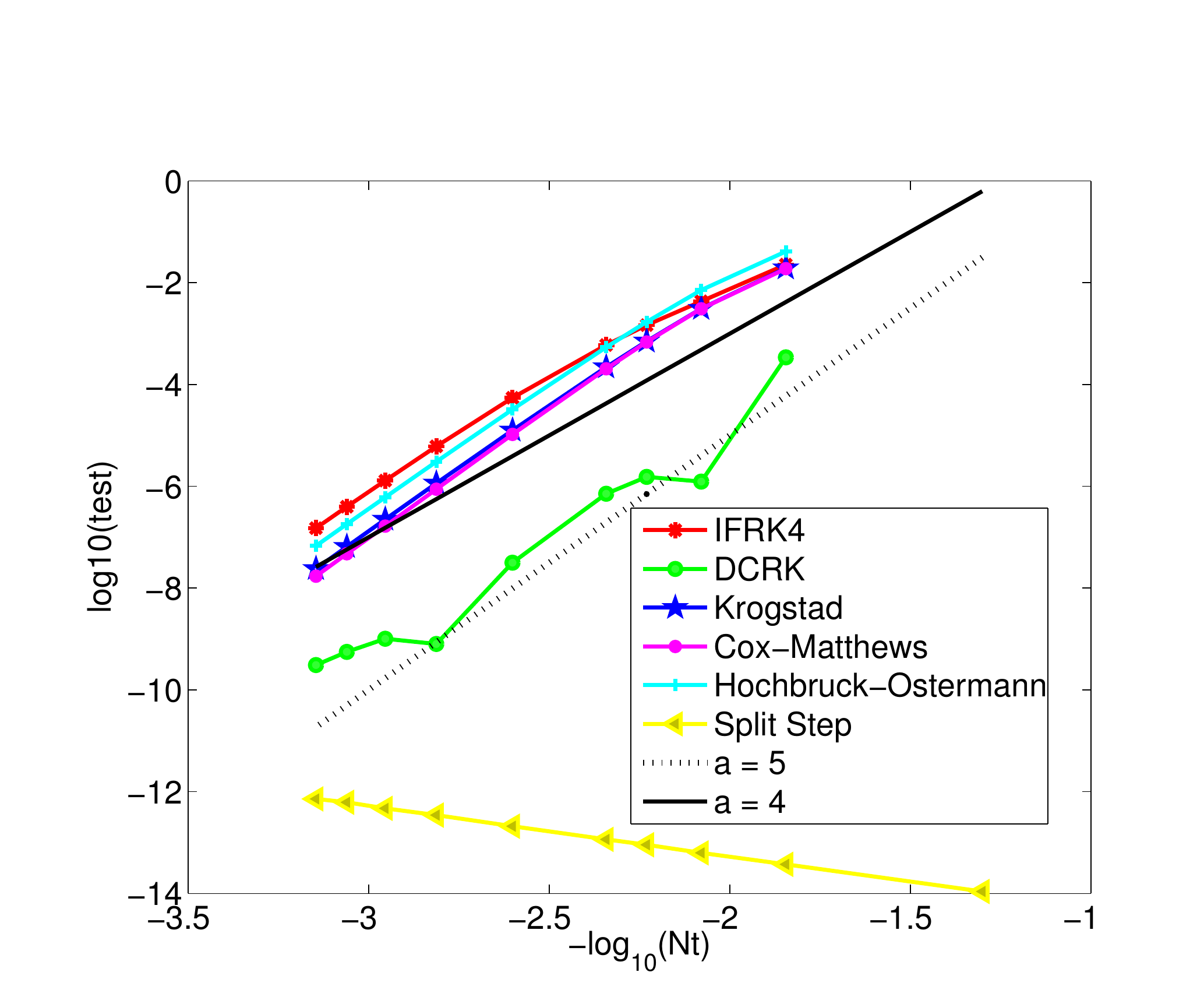}
 \includegraphics[width=0.5\textwidth]{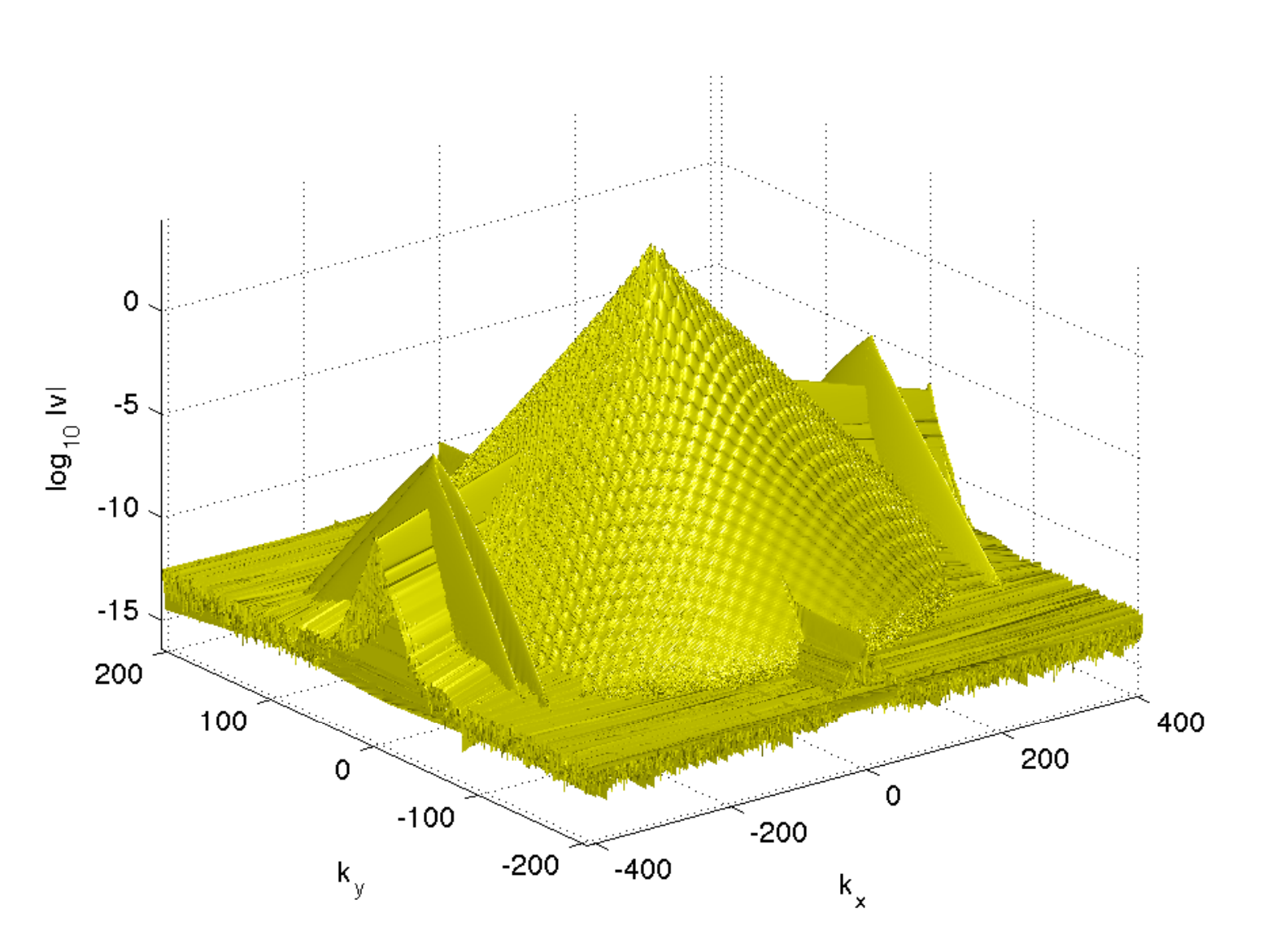}
 \caption{Non-conservation of the numerically computed $L_{2}$ norm 
 of the solution to the problem considered in 
 Fig.~\ref{figDSfocsmallreg} 
 in dependence on the time step
 (left) and the Fourier coefficients for the final time 
  (right). }
  \label{figDSfocsmallmass}
\end{figure}

For the small dispersion example for the focusing DS II equation of 
Fig.~\ref{figDSdefoc} it can be seen in Fig.~\ref{figDSfocsmallmass} 
that spatial resolution is almost achieved. There is a certain lack of 
resolution in the $k_{x}$ direction which leads to the formation of 
some structure close to $k_{y}=0$. This is related to the 
modulational instability of solutions to the focusing DS II equation. 
It will disappear for higher resolutions. Numerical conservation of 
the $L_{2}$ norm of the solution overestimates numerical accuracy by 
2-3 orders of magnitude for an error of the order of $10^{-3}$. Once 
more it cannot be used as an indicator for the numerical error in the splitting case, where it is almost 
independent of the time step. For the other cases numerical 
conservation of the $L_{2}$ norm
shows a dependence on $N_{t}$ between fourth and fifth order. This 
indicates as for the NLS case in \cite{ckkdvnls} that the numerical 
error has a divergence structure which leads to a higher order 
decrease of the $L_{2}$ norm than for the actual error. This behavior 
is also present in the above examples, but less pronounced. 

\section{Conclusion}
It was shown in this paper that fourth order time stepping schemes can be 
efficiently used for higher dimensional generalizations of the KdV 
and the NLS equations, where the stiffness of the system of ODEs 
obtained after spatial discretization can be a problem. Implicit 
schemes as IRK4 are computationally too expensive in the stiff regime, whereas standard explicit 
schemes as RK require for stability reasons too restrictive 
requirements on 
the time steps for the KP and DS equations. For these equations the 
non-localities in the PDEs lead to singular Fourier multipliers which 
make standard explicit schemes in practice unusable for stability 
reasons.

IMEX schemes do not converge in general for similar reason. 
Driscoll's composite RK variant is generally very efficient if the studied 
system is not too stiff, but fails to converge for strong stiffness. 
Exponential integrators do not have this problem.  The order reduction phenomenon is a 
considerable problem for IF schemes in the stiff regime, but less so 
for ETD schemes. The Hochbruck-Ostermann method performs in general 
best, but the additional stage it requires is in practice not 
worth the effort in comparison with Krogstad's or Cox-Matthews' method. The computation of 
the $\phi$-functions in ETD is inexpensive for the studied problems since it 
has to be done only once. 

Since stiffness is not the limiting factor for DS II, all schemes perform 
well in this context. But the modulational instability of the focusing case requires 
high spatial resolution we could not achieve in Matlab on the used 
computers for more general initial data. Thus the code will be 
parallelized to allow the use of higher spatial resolution without 
allocating too much memory per processor. 

\bibliographystyle{siam}
\bibliography{biblio}{}

\end{document}